\documentclass[11pt,a4paper,notitlepage,oneside]{amsart}
\usepackage{etoolbox}
\usepackage{microtype}
\usepackage{amsmath}
\usepackage{amsthm}
\usepackage[psamsfonts]{amssymb}
\usepackage[T1]{fontenc}
\usepackage{lmodern}
\usepackage[english]{babel}
\usepackage[utf8]{inputenc}
\usepackage{textcomp}
\usepackage{newunicodechar}
\newunicodechar{τ}{\ensuremath{\tau}}
\newunicodechar{ℝ}{\ensuremath{\R}}
\usepackage{enumitem}
\usepackage{mathtools}
\usepackage{caption,subfig}\usepackage[normalem]{ulem}
\usepackage{xcolor,soul}
\usepackage{scalerel}
\usepackage{relsize}
\usepackage[normalem]{ulem}
\usepackage{booktabs}

\usepackage[autostyle,italian=guillemets]{csquotes}
\usepackage[style=alphabetic,firstinits=true,maxbibnames=99,backend=biber]{biblatex}
\allowdisplaybreaks 
\usepackage{pgf,tikz,pgfplots}

\usepackage{mathrsfs}
\usetikzlibrary{arrows}

\usepackage{fancyhdr}
\usepackage{tikz-cd}
\usepackage{amscd}
\usepackage[a4paper, bottom=3cm, left=3cm,right=3cm]{geometry}
\usepackage{hyperref}
\hypersetup{hidelinks}
	
\usepackage{subfiles}
\usepackage{PhDefs}

\usepackage{mathtools}
\usepackage{dsfont}
\numberwithin{equation}{section} \theoremstyle{plain}

\newtheorem{thm}{Theorem}
\numberwithin{thm}{section}

\newtheorem{lemma}[thm]{Lemma}
\newtheorem{prop}[thm]{Proposition}
\newtheorem{cor}[thm]{Corollary}

\theoremstyle{definition}
\newtheorem{de}[thm]{Definition}

\newtheorem{rem}[thm]{Remark}

\theoremstyle{remark}
\newcounter{exc}

\numberwithin{exc}{subsection}
\newtheorem{exa}[thm]{Example}

\makeatletter
\renewenvironment{proof}[1][\textit{Proof}]{\par    \pushQED{\qed}  \normalfont \topsep6\p@\@plus6\p@\relax
  \trivlist
  \item[\hskip\labelsep
        \upshape
    #1\@addpunct{:}]\ignorespaces
}{  \popQED\endtrivlist\@endpefalse
}
\makeatother

\usepackage[disable]{todonotes}

\newcommand{\todoexp}[1]{\todo[color=blue!60!white!70!green!80!white,size=\footnotesize]{explain: #1}}

\newcommand{\todocre}[1]{\todo[color=purple!60,size=\footnotesize]{create: #1}}

\author{Ingmar Metzler}
\address{
	Department of Mathematics, ETH Zurich, 8092 Zurich, Switzerland
}
\email{ingmar.metzler@math.ethz.ch}

\title{A new converse theorem for Borcherds products}

\begin{document}
		\begin{abstract}
		We establish a new converse theorem for Borcherds products. 
		Moreover, the injectivity of the Kudla--Millson theta lift is demonstrated in the O(n,2) 
		case in greater generality than is currently available in the literature. 
		Both results are derived under the assumption of a single hyperbolic split of the base lattice.
	\end{abstract}
	\maketitle
		\tableofcontents
		\section{Introduction}

	In 1998, Richard Borcherds was awarded the Fields Medal in recognition of his pioneering contributions 
to the field of automorphic forms and mathematical physics \autocite{Borcherdsmathhistory}. 
He had introduced vertex operator algebras as well as a multiplicative lift of automorphic forms 
and employed these tools to prove the Conway--Norton moonshine conjecture. 
The constructed lift -- nowadays referred to as \emph{Borcherds lift} -- maps 
vector-valued elliptic modular forms to orthogonal automorphic forms with infinite product expansions. 
These \emph{Borcherds products} provide a diverse array of examples of orthogonal modular forms. 

To be more explicit, let $(L,\qq)$ be an even lattice of signature $(\possig,2)$ 
and $\Grassm$ be the associated Grassmannian which carries a natural complex structure. 
The natural arithmetic subgroup to consider is the \emph{discriminant kernel} $\disckernel[L] \leq \OG(L)$ 
which is a finite index subgroup stabilising the discriminant group of $L$.
On the elliptic side, let $\Df = L'/L$ and consider holomorphic modular functions 
$f\colon \Ha \to \C[\Df]$ of weight $k \in \Z/2$, 
transforming with the dual Weil representation $\overline{\rho}_{L}$ which may possess poles at the cusp. 
This space is referred to as the space of 
\emph{weakly holomorphic modular forms} and is denoted by $\MF_{L^{-},k}^{!}$. 
\\
Such a form $f \in \MF_{L^{-},k}^{!}$ admits a Fourier expansion 
\begin{equation*}	f(\tau) = \sum_{\mu \in L'/L} \sum_{m \in \Z - \qq(\mu)} a(\mu,m) e^{2 \pi i \tau m} \efr_{\mu}
\end{equation*}
where $\efr_{\mu}$ denotes the standard basis vector of $\C[\Df]$ being $1$ in the $\mu$ component and $0$ elsewhere. 
Assume $f$ has integral principal part, 
meaning that $a(\mu,m) \in \Z$ for all $m < 0$. 
Then the Borcherds lift associates to $f$ a meromorphic modular form $\Psi(z,f)$ on $\Grassm$ 
for $\disckernel[L]$ with unitary multiplier system of finite order such that 

\begin{enumerate}[label=\alph*)]
	\item 
	the weight of $\Psi(\argdot,f)$ is given by $a(0,0)/2$,
	
	\item 
	the divisor of $\Psi( \argdot , f)$ equals 
	\begin{equation*}		\frac{1}{2} \sum_{\mu \in L'/L}\sum_{n > 0} a(\mu,-n) Z(\mu,n),
	\end{equation*}
	where 
	\[
	Z(\mu,n) = \sum_{\substack{\lambda \in L + \mu \\ \qq(\lambda) = n}} \{z \in \Grassm \mid z \perp \lambda\}
	\]
	is the \emph{Heegner divisor}\index{Heegner!divisor} of discriminant $(\mu, n)$ (cf.\ Definition \ref{de:Heegner_divisor}), 	
	\item 
	the target form $\Psi(\argdot,f)$ has explicit infinite product expansions. 
\end{enumerate} 

As previously stated, 
this lift gives rise to a variety of new orthogonal modular forms and also to the monster denominator formula, 
which is utilised to prove the moonshine conjecture. 
In \cite[Pr.~16.10]{Borcherds1998}, Borcherds himself poses the question whether his lift may be reversed. 
That is, given a meromorphic modular form $F$ for $\disckernel[L]$ whose zeros and poles are supported on special divisors $Z(\mu,n)$, 
is there a form $f \in \MF_{L^{-},k}^{!}$ such that $\Psi(z,f)$ equals $F$ up to a non-zero constant factor? 
In this paper we provide the following answer (cf.\ Theorem~\ref{thm:converse_theorem_Borcherds_products}). 

\begin{thm}\label{thm:Intro_Conversetheorem}
	Assume that $L$ splits a hyperbolic plane and $\possig > 3$. 
	Then every meromorphic modular form $F$ with respect to $\disckernel[L]$ whose divisor is a linear combination of special divisors $Z(\mu,n)$ is (up to a non-zero constant factor) the Borcherds lift~$\Psi(z,f)$ 
	of a weakly holomorphic modular form $f \in \MF_{L^-,1-\possig/2}^{!}$. 
\end{thm} 

\noindent
The case $\possig = 3$ is treated separately in Section \ref{sec:KMlift_Injectivity}. 
More specifically, the same result holds when the Witt rank of~$L$ equals~$1$, 
while a weaker version is established for Witt rank~$2$. 
The result stated above improves upon the previously strongest known result \cite[Thm.~1.2]{BruinierConverse}, 
which either required an additional scaled hyperbolic split or needed a transition to a sublattice of \(L\). 
\\
The proof of Theorem~\ref{thm:Intro_Conversetheorem} is intricate. 
As a first step, we invoke a celebrated duality statement by Bruinier \cite{Bruinier2002}  
rendering the above converse statement equivalent to the injectivity of another theta lift with wide-ranging applications:  
the \emph{Kudla--Millson} lift.

\subsection*{The Kudla--Millson lift}

In the 1980s, Kudla and Millson \cite{Kudla1986} introduced special Schwartz forms $\KMform$ on the symmetric spaces attached 
to the classical groups $\OG(\possig,\negsig)$, $\UG(\possig,\negsig)$, and $\Sp(\possig,\negsig)$, 
taking values in closed differential forms. 
Their principal objective was to investigate cohomology classes of special cycles  by means of 
a theta correspondence, generalising the celebrated work of Hirzebruch and Zagier \cite{Hirzebruch1976} 
on Hilbert modular surfaces. 
More precisely, in the orthogonal case, Kudla and Millson symmetrised the Schwartz forms $\KMform$ 
over a base lattice $L$ 
and let the Weil representation of the symplectic group act on $\KMform$ in order to obtain a kernel
(cf.\ Definition~\ref{de:KMTheta})
\[
	\KMTheta(\tau,z), \quad \tau \in \Ha, \quad z \in \Grassm 
\]
in a symplectic and an orthogonal variable, which transforms automorphically in both variables. 
These may then be employed as an integral kernel to shift automorphic objects from the symplectic to the orthogonal setting or vice versa. 
Such an association is referred to as a \emph{theta lift}. 
The authors extended their work in \cite{Kudla1987} in great generality; for an introduction, see 
\cite{Kudla1990} and \cite{Bruinier2004a}. 
The applications of the Kudla--Millson theta correspondence range from the study of the cohomology 
of orthogonal and unitary Shimura varieties \cite{Kudla1986} and Arakelov theory of Shimura varieties \cite{Kudla2004a} over specific counting problems \cite{Engel2023},  
constructing mock modular forms and higher dimensional error functions \cite{Funke2017} to proving a converse theorem for Borcherds products \cite{Bruinier2002, BruinierConverse}. 

In this paper, we consider the Kudla--Millson lift on the space of vector-valued elliptic modular forms associated to the Weil representation $\dwrep_{L}$ of an even lattice $L$ of signature $(\possig,2)$. 
To be more explicit, let $k = \rk(L)/2$ and $\CF_{L,k}$ denote the space of cusp forms 
for the Weil representation $\dwrep_{L}$. 
Then  
\[
	\KMlift\colon \CF_{L,k} \to \Hforms^{(1,1)}(Y_{L}), \quad 
	f \mapsto \int_{\Mp_{2}(\Z) \backslash \Ha} \langle f , \KMTheta \rangle \Im^{k} \di \mu
\]
defines a linear map to the space of square-integrable harmonic differential forms of Hodge type~$(1{,}1)$ on the variety $Y_{L} = \disckernel[L] \backslash \Grassm$ and is referred to as the \emph{Kudla--Millson lift}. 
Here $\mu$ denotes the hyperbolic measure on the upper half-plane $\Ha$ on which $\Mp_{2}(\Z)$ acts via Möbius transforms. 
The question of its injectivity already arose in \autocite{Kudla1990} 
and may be used to compute the rational Picard number of the underlying Shimura variety \autocite{Bergeron2016}, 
as well as to derive properties of cones generated by special cycles \autocite{Bruinier2019} \autocite{Zuffetti2022}.  
However, 
the primary application we have in mind is the converse Theorem \ref{thm:Intro_Conversetheorem} for Borcherds products. 

There have been multiple results on the injectivity of the Kudla--Millson lift over the past two decades 
presented in \autocite[Thm.~5.12]{Bruinier2002}, \autocite[Cor.~4.11]{BruinierFunkesurj2010}, \autocite[Thm.~5.3]{BruinierConverse}, \autocite[Thm.~6.1]{Zuffetti2024}, and \autocite[Cor.~7.8]{Stein2023}. 
These advances are based on three fundamentally different methods. 
While there has been a recent success in unifying the results above by Zuffetti and the author \autocite[Thm.~6.2]{Metzler2026} by generalising Zuffetti's method, this has not shed any new light on the converse theorem for Borcherds products. 
The most potent result 
for this application remained Theorem \ref{thm:Bruinier_KMinj_conversepaper} by Bruinier 
which provides injectivity under the assumption that the lattice $L$ splits a hyperbolic plane 
as well as a scaled hyperbolic plane. 

In the context of this paper, 
we present a new approach to proving the injectivity of the Kudla--Millson lift by means of 
computing cycle integrals on the orthogonal variety. 
This approach eliminates the need for an additional scaled hyperbolic split, leading to the following result 
(cf.\ Theorem~\ref{thm:mainhyperbolicsplit}).

\begin{thm}\label{thm:Intro_KMinj}
	Let $(L,\qq)$ be an even lattice of signature $(\possig,2)$ with $\possig > 3$. 
	Assume that $L$ splits a hyperbolic plane. 
	Then the Kudla--Millson lift ${\KMlift}\colon \CF_{L,1 + \possig/2} \to \mathcal{H}^{(1,1)}(Y_{L})$ is injective.
\end{thm}

The case $\possig = 3$ is also discussed in Subsection \ref{ssec:Injectivity_hyp_split} and proven for the lattice~$L$ having Witt \mbox{rank $1$}. 
It should be stressed that in the case of no hyperbolic split, Bruinier constructed a subspace of the kernel 
and demonstrated that it is in general non-trivial even in 
the simple instance of lattices of prime level \cite[Sec.~6.1]{BruinierConverse}, disproving injectivity.

Our new proof strategy is to compute special cycle integrals of Kudla--Millson liftings~$\KMlift(f)$ 
of cusp forms~$f \in \CF_{L,k}$. 
For suitably chosen cycles, these may be expressed in terms of special values of certain $L$-series $L(f,s)$
associated to~$f$. 
If~$\KMlift(f)$ is assumed to vanish, its cycle integral must vanish as well, 
yielding zeroes of the aforementioned $L$-series. 
As a consequence, relations between Fourier coefficients of $f$ appearing in the construction of $L(f,s)$ are established. 
By constructing an infinite family of subseries of these special \(L\)-values -- via the same cycle integrals -- 
we isolate the coefficients, ultimately showing that each of them must vanish.  
Since the Kudla--Millson lift is linear, injectivity is inferred. 
As a consequence, we deduce a converse theorem for Borcherds products under minimal splitting assumptions. 

\subsection*{Acknowledgments}
I am grateful to my supervisor Jan Bruinier for his guidance and numerous fruitful discussions exceeding the scope of this project. 
Also, I would like to thank Markus Schwagenscheidt and Paul Kiefer for substantial improvements to the manuscript. 
Further, I would like to extend my gratitude to Özlem Imamoglu and Jens Funke for helpful remarks.
The author was supported by the CRC TRR~$326$ `GAUS', project number~$444845124$, the LOEWE project `USAG', and ETH Zurich.

		\section{Vector-valued modular forms}

We briefly present the setting for this paper, including lattices, 
discriminant forms, the Weil representation, and elliptic modular forms. 

\subsection{The classical setting} 
Let $(L,\qq)$ denote a non-degenerate even quadratic lattice of signature $(\possig,\negsig)$ with associated bilinear form $\bilf$ such that $\qq(x) \coloneq \bilf(x,x)/2$. 
$\Hom_{\Z}(L,\Z)$ may be identified via $\bilf$ with a lattice $L' \subset L \otimes_{\Z} \Q$ referred to as the \emph{dual} and we write $\level(L)$ 
for the level of the lattice, being the smallest positive integer such that $\bilf(L',\level(L) \cdot L')$ 
is even. Recall the inclusion $L \leq L'$ and note that $\Df \coloneq L'/L$ is a finite abelian group 
whose order is the absolute value of the determinant of a representing matrix of $\bilf$. 
Further, the quadratic form $\qq$ descends to a quadratic form $\qd\colon \Df \to \Q/\Z$ that is non-degenerate 
and such a pair $(\Df,\qd)$ is called a \emph{discriminant form}. 
In fact, any discriminant form arises from a non-degenerate even lattice in the above fashion. 
If we set $\sig(\Df) = \possig - \negsig \mod 8$ for some lattice $L$ inducing $\Df$, 
this quantity is well defined by Milgram's formula. 
For a natural number $n \in \N$ denote by $L(n)$ the pair $(L, n \cdot \qq)$ and set~$\Df(n) \coloneq L(n)'/L(n)$. 

Additionally, write $\Dftors{n}$ for the $n$-torsion in $\Df$ and $\Dfmult{n}$ for the $n$-multiples. 
We have the exact sequence 
\begin{equation}\label{eq:Df_torsion_mults}	
	\begin{tikzcd}
		0 \arrow[]{r}{} 		& \Dftors{n} \arrow[]{r}{\iota}		&	 \Df \arrow[]{r}{\cdot n}		& 	\Dfmult{n}  \arrow[]{r}{}		&	0 
	\end{tikzcd}
\end{equation}	
and find that $\Dfmult{n}$ is exactly the orthogonal complement of $\Dftors{n}$. 
We denote by $\C[\Df]$ the group algebra of $\Df$ with standard basis $\efr_{\lambda}$ for $\lambda \in \Df$.
There is a representation of $\Gamma(1) \coloneq \SL_{2}(\Z)$ on the group algebra $\C[\Df]$ first described by Schoeneberg. 
In fact, to obtain a well-defined representation when $\sig(\Df)$ is odd, 
it must be considered on the metaplectic extension, rather than on $\SL_{2}$ directly. 
Let $\MGLtwop(\R)$ denote the metaplectic double cover of $\GL_{2}^{+}(\R)$. 
A model is given by pairs $(\gamma,\phi)$, where $\gamma \in \GL_{2}^{+}(\R)$ and $\phi\colon \Ha \to \C$ 
is a holomorphic square root of the factor of automorphy $j(\gamma,\tau)$, meaning $\phi^2(\tau) = j(\gamma,\tau)$. 
Multiplication is then declared by 
\[
	(\gamma, \phi(\tau)) (\gamma',\phi'(\tau)) = (\gamma \gamma', \phi(\gamma' \tau) \phi'(\tau)). 
\]
Projection to the first component \( \MGLtwop \ni (\gamma, \phi) \mapsto \gamma \in \GL_{2}^{+}(\R)\) defines a covering. 
The action of elements $\abcds = \gamma \in \GL_{2}(\R)$ on the upper half plane $\Ha$ is given by Möbius transforms 
\(\tau \mapsto (a \tau + b)/(c \tau + d)\) with diagonal invariance group $\R^{\times} \into \GL_{2}^{+}(\R)$ and induces an action on $\MGLtwop(\R)$ via pullback. 
For the preimage of $\SL_{2}(\Z)$ in $\MGLtwop(\R)$, we write $\Mp_{2}(\Z)$. 
This group is generated by the two elements $\Telmp = \left( \sm{1}{1}{0}{1}, 1\right)$ 
and $\Selmp = \left( \sm{0}{-1}{1}{0}, \sqrt{\tau} \right)$, where $\sqrt{\tau}$ denotes the standard holomorphic square root. 
These elements project to the standard generators of $\SL_{2}(\Z)$ and $\Selmp^2 = \Zelmp \coloneq \left(- \sm{1}{0}{0}{1} , i \right)$ 
generates the centre of $\Mp_{2}(\Z)$. 
\\
We now turn towards describing the Weil representation.
It suffices to describe the action of the Weil representation $\rho_L$ associated to $L$ on the generators $\overline{T}$ and $\overline{S}$ of $\Mp_{2}(\Z)$. 
	We abbreviate $e(\argdot) \coloneq \exp(2 \pi i \argdot)$ and note that these operations are given by 
\begin{align}
	\dwrep_L(\Telmp) \mathfrak{e}_\lambda 	&=	e(\qd(\lambda)) \cdot \mathfrak{e}_\lambda,	\label{eq:op_T_via_rho}\\
	\dwrep_L(\Selmp) \mathfrak{e}_\lambda 	&=	\frac{e(-\sig(\Df)/8)}{\sqrt{|\Df|}} \cdot \sum_{\mu \in \Df} e(-\beta(\lambda,\mu)) \mathfrak{e}_\mu \label{eq:op_S_via_rho}.
\end{align}
It is apparent that the representation $\rho_{L}$ only depends on the discriminant form $\Df$.
\begin{de}\label{de:dwrep}
	Let $(\Df,\qd)$ be a discriminant form. 
	Then $\dwrep_{\Df} \coloneq \dwrep_L$ is a unitary representation of $\Mp_{2}(\Z)$ on $\C[\Df]$ that is fully determined by 
	\eqref{eq:op_T_via_rho} and \eqref{eq:op_S_via_rho}. 
\end{de}
Note that the generator of the centre $\overline{Z}$ operates for $\lambda \in \Df$ as 
		\begin{align}
			\dwrep_\Df(\overline{Z}) \mathfrak{e}_\lambda = e(-\sig(\Df)/4) \cdot \mathfrak{e}_{-\lambda}. \label{eq:op_Z_via_rho}
		\end{align}
		As a consequence, the Weil representation $\rho_\Df$ factors through $\Mp_{2}(\Z) / \langle\overline{Z}^2\rangle \simeq \SL_2(\Z)$, 
		if, and only if, $\sig(\Df) \equiv 0 \mod 2$.
		Further, writing $N \coloneq \level(\Df)$ we find that $\dwrep_{\Df}$ is trivial on the congruence subgroup $\Gamma(N)$ for even rank, 
	so that it factors through the finite group
	\[
		\SL_{2}(\Z /N\Z) \simeq \Gamma(1) / \Gamma(N).
	\]
	In case of odd rank, the oddity formula \cite[Chap.~15 p.~383 (30)]{ConwaySloane1998} implies $4 \mid N$, in particular, $\Df$ contains $2$-adic Jordan components. 
	In this case, there is a well known section 	\begin{align}\label{eq:section_Gamma4_metaplectic}
		s\colon \Gamma(4) \to \Mp_2(\Z), \qquad \gamma = \abcd \mapsto \left( \abcd , \epsilon_{d}^{-1} \left( \frac{c}{d} \right) \sqrt{j(\gamma, \tau)} \right)
	\end{align}
	with $\epsilon_{d} = 1$ or $i$, depending on whether $d \equiv 1$ or $3 \mod 4$. 
		The argument at the end of \cite[Thm.~5.4]{Borcherds2000} yields that $\dwrep_{\Df}$ is trivial on $s(\Gamma(N))$ and factors through the central extension of $\SL_{2}(\Z/N\Z)$ by $\{\pm 1\}$ given by
	\[
		\Mp_{2}(\Z) / s(\Gamma(N)). 
	\]

The natural projection $\Mp_{2}(\Z) \to \SL_{2}(\Z)$ has the section $\gamma \mapsto \tilde{\gamma} \coloneq (\gamma, \sqrt{j(\gamma,\tau)})$, 
where we selected the standard branch of the holomorphic square root. Note that this mapping is not a group homomorphism. 
Before continuing we briefly recall the notion of a classical modular form. 
\begin{de}\label{de:Peterssonslash}
	Let $\gamma \in \GL_{2}^{+}(\R)$ and $f\colon \Ha \to \C$ be a function. Define 
	\[
		f\vert_{k} \gamma (\tau) \coloneq \det(\gamma)^{k/2} j(\gamma, \tau)^{-k} f(\gamma \tau).
	\]
\end{de}

Finally, a holomorphic modular form $f\colon \Ha \to \C$ of weight $k$ for some congruence subgroup $\Gamma \leq \Gamma(1)$ with character $\chi$ on $\Gamma$  
is a holomorphic function such that for all $\gamma \in \Gamma$ we have $f \vert_{k} \gamma = \chi(\gamma) f$ and $f$ is holomorphic at the cusps of $\Gamma$. 
The space of such functions is denoted by $\MF_{k}(\Gamma,\chi)$ and the subspace of elements vanishing at the cusps is denoted by $\CF_{k}(\Gamma,\chi)$ 
and called the space of cusp forms. 
For the technical modifications required for forms of half-integral weight, we refer to \cite[p. 444]{ShimuraMF_HalfIntegralWeight}.

\begin{de}
	A function $f\colon \Ha \to \C[L'/L]$ is called \emph{(vector-valued) modular function} of weight $k \in \frac{1}{2} \Z$ if for all $\gamma = (M, \phi) \in \Mp_{2}(\Z)$ the following transformation law holds 
	\[
	f(\gamma \tau) = \phi(\tau)^{2k} \cdot \left[\rho_L(\gamma) f\right](\tau).
	\]
\end{de}

The operation of $\Telmp$ implies that such a function $f$ possesses a Fourier expansion if it is holomorphic,
meaning there are coefficients $a(\lambda,n) \in \C$ such that 
\begin{equation}\label{eq:VVMF_FourierExp}
	f( \tau ) = \sum_{\lambda \in L'/L} \sum_{n \in \qd(\lambda) + \Z} a(\lambda,n) e(n\tau) \efr_{\lambda}.
\end{equation}
In case the discriminant group $L'/L$ is trivial, the index $\lambda \in L'/L$ is eliminated from the notation. 
We say that a holomorphic modular function $f$ is a \emph{modular form}, 
if all coefficients $a(\lambda,n)$ in \eqref{eq:VVMF_FourierExp} vanish for negative $n \in \Q$.  
Further, it is called a \emph{cusp form}, if all coefficients $a(\lambda,n)$ vanish for non-positive $n \in \Q$. 
The spaces of such modular functions are denoted by $\MF_{L,k}$ and $\CF_{L,k}$. 
Given a modular form $f \in \MF_{L,k}$, there is a general procedure to induce a modular form to a sublattice $M \leq L$. 
In that case, $M \leq L \leq L' \leq M'$ implying the inclusion $L'/M \leq M'/M$ and 
there is a lift $\uplift{L}{M}\colon \MF_{L,k} \to \MF_{M,k}$ declared by 
\begin{align}\label{eq:upop_concrete}
    \left(\uplift{L}{M}f\right)_\mu = 
    \begin{cases}
        f_{\overline{\mu}}, 	& \textnormal{if } \mu \in L'/M, \\ 
        0,						&	\textnormal{otherwise}, 
    \end{cases}
\end{align} 
where $\overline{\mu}$ refers to the image of $\mu \in L'/M$ under the projection to $L'/L$ (cf. \cite[Sec.~4]{Scheithauer2015}). 
By computing its adjoint, we obtain an operator \(\downlift{M}{L}\colon \MF_{M,k} \to \MF_{L,k}\) with 
\begin{align}\label{eq:downop_concrete}
(\downlift{M}{L}g)_{\overline{\mu}} = \sum_{\nu \in L/M} g_{\mu + \nu} 
\end{align}
where $\mu$ denotes an arbitrary but fixed lift of some $\overline{\mu} \in L'/L$ to $M'/M$. 

There are also more concrete examples of modular forms, theta functions, in particular. 

\begin{de}\label{de:Theta}
	Let $L$ be even positive definite, $\lambda \in L'/L$, and $\tau = u + iv \in \Ha$. 
	Define 
	\begin{align}
			\theta_{L,\lambda}(\tau)
			\coloneq\;
			& \sum_{l \in L + \lambda} e\left(\tau \qq(l) \right),  
			\qquad 
			\Theta_{L} \coloneq \sum_{\lambda \in L'/L} \theta_{L,\lambda}(\tau) \efr_{\lambda}.
	\end{align}
\end{de}

These theta series $\Theta_{L}$ are modular forms of weight $\possig/2$ (cf.\ \cite[Thm.~4.1]{Borcherds1998}) 
and their Fourier coefficients equal representation numbers of the lattice $L$. 

\begin{rem}\label{rem:downarrow_translates_thetafunctions}
	The operator $\downlift{M}{L}$ 
	translates the theta function associated with $M$ to a theta function of $L$. 	In the notation of Definition~\ref{de:Thetafun_Schwartzform} this means 
	\begin{align*}
		\downlift{M}{L}\left(\Theta_{M}(\tau,z;\varphi) \right) 	
		&=\; 	\downlift{M}{L}\left(\sum_{ \lambda \in M'/M} \sum_{x \in M + \lambda} (\omega_\infty (g_\tau) \varphi)(x,z) \e_\lambda\right) 	\\
		&=\;	\sum_{ \lambda \in L'/L} \sum_{x \in L + \lambda} (\omega_\infty (g_\tau) \varphi)(x,z) \e_\lambda 	\\
		&=\;	\Theta_{L}(\tau,z;\varphi).	
	\end{align*}
\end{rem}

\begin{prop}\label{prop:arrowopsrelation}
		For $M \leq L$ even non-degenerate lattices and $k \in \tfrac{1}{2}\Z$, there are linear operators 
	\[
	\downlift{M}{L}\colon \MF_{M,k} \to \MF_{L,k}, \qquad \uplift{L}{M}\colon \MF_{L,k} \to \MF_{M,k}
	\]
	fulfilling the following properties. 
	\begin{enumerate}[label=\alph*)]
		\item 
		The operators map modular (cusp) forms to modular (cusp) forms.
		
		\item 
		The lift $\downlift{M}{L}$ is compatible with theta functions, i.e.\ $\downlift{M}{L}(\Theta_{M}) = \Theta_{L}$.
		
		\item 
		The operators are adjoint to each other on the space of cusp forms:
		\[
		\langle \, \cdot \, ,  \downlift{M}{L}(\, \cdot \,) \rangle_L = \langle \uplift{L}{M}(\, \cdot \,) , \, \cdot \, \rangle_M. 
		\]
		This identity also holds, if one of the cusp forms is replaced by a modular form.
		
		\item 
		The image $B \coloneq \im(\uplift{L}{M})$ is the subspace consisting of functions that are supported on $L'/M$.
		We find  
		\begin{align*}
			\downlift{M}{L} \circ \uplift{L}{M} 				&= |L/M| \cdot \id, \\
			\uplift{L}{M} \circ \downlift{M}{L} \big\vert_{B} 	&= |L/M| \cdot \id_B. 
		\end{align*} 
	\end{enumerate}
	\end{prop}

Beyond the theta series presented above, Eisenstein series represent
another significant class of modular forms. We refer to \cite[2.5 pp.~30-36]{Kiefer2021} and {\cite{BruinierKuehn2003}}.

\begin{de}\label{de:VVEisensteinseriesnonholomorphic}
	Let $\lambda \in \Df$ be isotropic and $k \in \Z/2$. 
	Define the \emph{Eisenstein series}\index{Eisenstein series!non-holomorphic} 
	\begin{equation}
		E_{L,\lambda,k}(\tau,s) \coloneq \frac{1}{2} \sum_{\gamma \in \overline{\Gamma_{\infty}} \setminus \Mp_{2}(\Z)} \left( \Im(\tau)^s \efr_\lambda \right) \vert_{L,k} \gamma .
	\end{equation}
\end{de}

By construction this function is $\Mp_{2}(\Z)$ modular of weight $k$.

\begin{thm}\label{thm:ELk_convergence_hypLaplacian_Eigenfunction}
	The Eisenstein series $E_{L,\lambda,k}$ converges normally on $\Ha$ for $\Re(s) > 1 - \frac{k}{2}$, 
	is real analytic in $\tau$ and an eigenfunction of the hyperbolic Laplace operator of weight $k$ with eigenvalue $s(s + k -1)$. 
\end{thm}

\begin{rem}\label{rem:VVEisenstein_meromorphiccont}
	Since~$\dwrep_{L}$ acts trivially on $\Gamma(\level(L))$, 
	each component of $E_{L,\lambda,k}(\tau,s)$ is a scalar-valued Eisenstein series for $\Gamma(\level(L))$,  
	which admits meromorphic continuation in $s$. 
\end{rem}

We refer to \cite{Metzler2025} for the following comment on asymptotic bounds required for convergence at a later stage. 
\todo{I am unhappy with this, as it does not suffice for the convergence statement which is required later. But the other one is complicated.}

\begin{rem}\label{rem:Eis_locunifbound}
	Let $s \in \C$ such that $E_{L,\lambda,k}(\tau,s)$ is holomorphic. 
	We find for $\tau = u + iv \in \Ha$ and $\sigma = \max\{\Re(s), \Re(1-k-s)\}$ that 
	\begin{equation}
		E_{L,\lambda,k}(\tau,s) = \LandO(v^{\sigma}), \qquad v \to \infty.
	\end{equation}
	In fact, the constant required to bound may be chosen locally uniformly in $s$.
\end{rem}

\subsection{Eisenstein series}\label{sec:ES}

There is also a representation-theoretic approach to Eisenstein series, described succinctly in \cite{KudlaYangSL2} which is suitable for applying the Siegel--Weil formula. 
We recall here the essentials of the construction, following their notation, 
and refer to the source for details. 

For the group $G_p=\SL_{2}(\Q_{p})$ the upper triangular matrices $P$ form a maximal parabolic subgroup, its standard Borel subgroup.  
We install the following notation: 
\begin{align}
	M &\coloneq \left\{ \m a 0 0 {a^{-1}} =: M(a) \; \middle\vert \; a \in \Q_{p}^{\times}  \right\}, \quad 
	N &\coloneq \left\{ \m 1 b 0 1 =: N(b) \; \middle\vert \; b \in \Q_{p}  \right\}.
\end{align}
Note that $\SL_{2}(\A)$ admits the so called \emph{Iwasawa} decomposition into \(N(\A) M(\A) K\) \index{Iwasawa!decomposition} 
with \(K = K_\infty K_f \), where \(K_\infty = \SO_2(\R)\) and \( K_f = \SL_2(\hat{\Z}) \).
Write 
\[
	k_\vartheta = \m{\cos(\vartheta)}{\sin(\vartheta)}{-\sin(\vartheta)}{\cos(\vartheta)} \in K_{\infty} \coloneq \SO_{2}(\R)
\] 
and \(\tau = u + iv \in \Ha\), then
\begin{align}
g_{\tau}k_{\vartheta} i 
= 
\tau 
= 
u + iv \in \Ha. \label{eq:Iwasawa_decomposition}
\end{align}
This decomposition lifts to \(\MSL(\A)\), the metaplectic extension of $\SL_2(\A)$ by the torus \(\Tc\), 
where the respective preimages of groups in \(\MSL(\A)\) are marked by an apostrophe and we assume the same parametrisation as in \cite{KudlaYangSL2}. 
Additionally, recall that there is a two-fold cover, denoted by \(\Mp_{2}(\A)\) sitting inside \(\MSL(\A)\). 
Moreover, note that for a rational quadratic space \((V,\qq)\), 
there is a Weil representation \(\omega_{p}\) of \(\MSL(\Q_{p})\) on the space of Schwartz--Bruhat functions \(\Sspace(V_p)\), where $V_p = V \otimes \Q_p$. 
Additionally, we will require the characters associated to the quadratic space induced by the Hilbert symbol \(\Hsym_{p}\): 
\begin{equation}\label{eq:Hilbertsymcharacter}
	\Q_p^{\times} \ni a_p \mapsto \Hsym_{p}(a_p , \disc(V_p)), \qquad \A^{\times} \ni (a_p) \mapsto \Hsym_{\A}(a, \disc(V)) \coloneq \prod_{p} \Hsym_{\Q_p}(a_p,\disc(V_p)). 
\end{equation}
For a more comprehensive description of the setting, we refer to \cite[Sec 8.5]{KRY2006}.

For an idele class character \(\chi \in \Q^\times \backslash \A^\times\) there are two principal series representations on \(\MSL(\A)\). 
These are the even and odd representation \(I(s,\chi)\) and their sections fulfil
\begin{align}
		\Phi(p'g',s) = \chi(a) \cdot |a|^{s+1} \Phi(g',s) 
		\begin{cases} 
			1,	&	\textnormal{in the even case,}	\\	
			z,	&	\textnormal{in the odd case,}
		\end{cases}\label{eq:princseriesrepsectiontransformation}
\end{align}
where $p' = [n(b)m(a),z] \in P_{\A}'$ is in the standard Borel subgroup of \(\MSL(\A)\)
and the same parametrisation as in \cite[p.~287]{KRY2006} is chosen. 	
A section $\Phi(s)$ is called \emph{standard}\index{section!standard}\index{standard!section} if its 
restriction to the maximal compact subgroup $K_{\A}<\SL_{2}(\A)$ ($K_{\A}'<\MSL_{2}$, respectively) 
is independent of $s$ and \emph{factorisable}\index{section!factorisable}\index{factorisable section} 
if $\Phi = \otimes_{p} \Phi_{p}$ is a primitive tensor for a family of sections \(\Phi_p \in I_p(s,\chi)\) in the associated local principal series representation. 
Note that every section $\Phi(g,s) \in I(s,\chi)$ is determined by its values on $K_{\A}'$ which is a consequence of the Iwasawa decomposition.

The simplest example of a section at a finite place is the normalised standard section. 

\begin{exa}\label{exa:sphericalsectionprincipalseriesp}
	Let $2 < p < \infty$ and note that there is a section of $K_p = \SL_{2}(\Z_p) \into \MSL_{2}(\Q_p)$.
	Then the standard section $\Phi_{p}^{0} \in I_p(s,\chi)$ determined by the condition \(\Phi_{p}^{0}(k,s) = 1\)
	for all $k \in K_p = \SL_{2}(\Z_p)$ is called the \emph{spherical section}\index{section!spherical}\index{spherical section}  
\end{exa}

Recall that the characters of the double cover $\MSO_{2}(\R)$ sitting in $K_{\infty}'$ are of the following form (cf.\ \cite[Sec.~2]{Bruinier2004a}). 

\begin{rem}\label{rem:MSO_characters}
	For every $l \in \Z/2$ there is a character of $\MSO_2(\R)$ denoted by $\nu_l$ with values  
	\[
		\nu_l\left( k_\vartheta , \pm \sqrt{j(k_\vartheta,\tau)} \right) \mapsto \pm \sqrt{j(k_\vartheta,i)}^{-1} = \pm e^{il \vartheta }.
	\]
\end{rem}

With this family of characters, one may define concrete elements of the principal series representation at the infinite place.  

\begin{exa}\label{exa:standardsecweightlprincipalseriesinfty}
	Assume $\chi$ is a quadratic character. 
	For $l \in \Z/2$ imposing the conditions 
	\[
	\Phi_{\infty}^l(gk,s) = \nu_l(k) \Phi_{\infty}^l(g,s), \quad \Phi_{\infty}^l(1,s) = 1
	\]
	together with~\eqref{eq:sectioninfty_l_paritycondition} yields the \emph{normalised eigenfunction}\index{normalised eigenfunction} 
	of weight $l$. These functions span the $K_\infty$ finite functions in $I_\infty(s,\chi)$. 
		
	For $\tau = u + i v \in \Ha $ we find 
	\begin{align*}
		\Phi_{\infty}^l(g_{\gamma \tau} k_{\vartheta},s) 
		&=	\Im(\tau)^{l/2}	 \frac{\Im(\gamma\tau)^{(s+1-l)/2}}{j(\gamma, \tau)^l}.
	\end{align*}
\end{exa}

There is an intertwining map which may be used to construct elements of the principal series representations 
by means of the Weil representation. 
In fact, when considering the action of the Weil representation   
and comparing it to \eqref{eq:princseriesrepsectiontransformation}, the following statement is inferred.

\begin{lemma}\label{lem:Intertwiningintoprincseries}
	Let $(V,\qq)$ be a quadratic $\Q$ space of dimension $\dimV$ with character $\chi_{V} = (\argdot, \disc(V))$.
	Further, let $p \leq \infty$, $\omega_p$ be the associated Weil representation and set 
	$s_0 \coloneq \frac{\dimV}{2} - 1$. 
	Then 
	\begin{align}\label{eq:lem:Intertwiningintoprincseries}
		\lambda_{p}\colon \Sspace(V_p) \to I_p(s_0,\chi), \qquad \lambda_{p}(\varphi)(g') \coloneq \omega_p(g')\varphi(0)
	\end{align} 
	is a well defined map. 
	In fact, it is almost always surjective, except for three special cases that are described in \cite[Sec.~4 p.~2286]{KudlaYangSL2}.
\end{lemma}

Clearly, to any of the sections arising from the above intertwining operator \eqref{eq:lem:Intertwiningintoprincseries}, there is an associated canonical standard section.
This method reproduces the already known examples from above but also yields new examples.

\begin{exa}\label{exa:Principalseriesec_from_intertwining}
	\begin{enumerate}[label=\alph*)]
		\item 
		Let $V_{\infty}$ have signature $(\possig,\negsig)$. 
		Then for any element $z$ of the associated Grassmannian, selecting the positive standard majorant\index{standard!majorant} 
		$\qq_{z}^{+}$ of $\qq$  (cf.\ \eqref{eq:standardmajorant}) and setting 
		\begin{equation}
			\varphi_{\infty}(x,z) \coloneq e^{- 2 \pi \qq_{z}^{+}(x)} \in \Sspace(V_{\infty})
		\end{equation}
		associates the standard section from Example~\ref{exa:standardsecweightlprincipalseriesinfty} 
		of weight $l = \tfrac{\possig-\negsig}{2}$ via the intertwining operator of Lemma~\ref{lem:Intertwiningintoprincseries} in the archimedean case.
				
		\item 
		Let $L_p$ be unimodular, then the standard section associated to $\lambda_{V_p}(\charfun_{L})$ equals
		the spherical section $\Phi_{p}^{0}(s)$ from Example~\ref{exa:sphericalsectionprincipalseriesp}.

	\end{enumerate}
\end{exa}

Further, we find an operator for the finite adelic part 
\begin{equation}\label{eq:intertwiningoperator_finitepart}
	\lambda_{f}\colon \Sspace(V_{\A_{f}}) \to I_{f}(s_0,\chi) = \otimes_{p < \infty}' I_{p}(s_0,\chi), \qquad \lambda_{f}(\varphi)(g_{f}') = \omega_{f}(g_{f}') \varphi(0).
\end{equation}
For $\mu \in L'/L \simeq \hat{L}'/\hat{L}$ there is the related element $\varphi_{\mu} \coloneq \mathds{1}_{\mu + \hat{L}} \in \Sspace(V_{\A_{f}})$ for which we find 
$\lambda_{f}(\varphi_{\mu}) \in I_{f}(s_0,\chi)$. 
The associated standard section will be denoted by $\Phi_{\mu}$ and plays a major role in 
the construction of Eisenstein series below. 

With the tools described above, we may turn towards Eisenstein series.

\begin{de}
	For a standard section $\Phi(g,s) \in I(s, \chi)$ define 
	the \emph{Eisenstein series}\index{Eisenstein series!adelic}\index{Adele!Eisenstein series}
	\begin{equation}\label{eq:EisensteinSeries}
		\Eisadel(g,s;\Phi) := \sum_{\gamma \in P(\Q) \setminus \SL_{2}(\Q)} \Phi(\gamma g, s).
	\end{equation}
\end{de}
Note that this series converges for $\Re(s) > 1$, is, by definition, left invariant under $\SL_{2}(\Q)$, and has a meromorphic continuation to the whole complex plane in $s$ 
due to Langlands (cf\ \cite[6 Lem~6.1]{Langlands2006}).

We are interested in concrete instances of the above series and restrict, from now on, 
to a quadratic character $\chi$ which may be realised by a square free integer $d \in \Z$ via 
\begin{align*}
	\chi(x) = 
	\begin{cases}
		\Hsym_{\A}(x,d), & \textnormal{even case}, \\
		\Hsym_{\A}(x,2d), & \textnormal{odd case}.
	\end{cases}
\end{align*}
Here, $\Hsym_{\A}$ denotes the Hilbert symbol and 
the even or odd case refers to the principal series representation factoring through $\SL_{2}$ or not.\footnote{This will agree with the signature being even or odd in our later setting.}  
Recall that in order for the section $\Phi_{\infty}^{l}$ of Example~\ref{exa:standardsecweightlprincipalseriesinfty} to be non-trivial and hence part of $I_{\infty}(s,\chi)$, 
the following parity condition has to be satisfied  
\begin{align}\label{eq:sectioninfty_l_paritycondition}
	\begin{cases}
		(-1)^{l} = \sig(d), 		&	\textnormal{even case,} \\
		l \equiv \sig(d)/2 \mod 2, 	&	\textnormal{odd case.}
	\end{cases} 
\end{align}

In this context, classical Eisenstein series are representable. 
First, recall that there is a section $P_{\infty} \to \Mp_{2}(\R)$ of the standard Borel subgroup, 
so that $g_{\tau}$ may directly be interpreted as an element in $\Mp_{2}(\R)$.

\begin{exa}
	For $\chi$ being the trivial character and $g_{\tau}$ as above we find by Example~\ref{exa:standardsecweightlprincipalseriesinfty} and Example~\ref{exa:sphericalsectionprincipalseriesp} 
	with the definition $\Phi_{f}^{0} = \otimes_{p} \Phi_{p}^{0}$ that 
	\[
		\Eisadel(g_{\tau},s;\Phi_{\infty}^l \otimes \Phi_{f}^{0}) = \Im(\tau)^{l/2} \sum_{\gamma \in \Gamma_{\infty} \backslash \SL_{2}(\Z)} \frac{\Im(\gamma\tau)^{(s+1-l)/2}}{j(\gamma,\tau)^l}.
	\]
\end{exa}

\todocre{Declare standard injection of $\SL_{2}$ into $\Mp_{2}$ (ofc not as groups) as we did with Riccardo}

The preceding example demonstrates how classical Eisenstein series arise from the adelic construction. 
This motivates considering the Eisenstein series in \eqref{eq:EisensteinSeries} 
as functions on the upper half plane. 
Recall that any factorisable finite standard section $\Phi_{f}(g,s) \in I_{f}(\chi,s)$ is invariant under some open subgroup $K_0$ of $\SL_{2}(\hat{\Z})$.
Hence, by strong approximation the following Eisenstein series determines the series $\Eisadel(g',s;\Phi_{\infty}^{l} \otimes \Phi_{f})$ completely. 

\begin{de}\label{de:Eisadel_ongroup_toEisadel_intau}
	For $\Phi_\infty^l$ from Example~\ref{exa:standardsecweightlprincipalseriesinfty} transforming with $\nu_l$ under $\overline{K_\infty}$ and a finite standard section $\Phi_{f} \in I_{f}(\chi,s)$, define
	\begin{equation}
		\Eisadel(\tau,s; \Phi_\infty^l \otimes \Phi_f) := \Im(\tau)^{-l/2} \cdot \Eisadel(g_{\tau},s;\Phi_\infty^l \otimes \Phi_f).
	\end{equation}
\end{de}

The series $\Eisadel(\tau,s; \Phi_\infty^l \otimes \Phi_f)$ then defines a real analytic modular form of 
weight~$l$ and we have seen how to realise classical scalar valued Eisenstein series. 
However, in this paper, vector-valued modular forms play a major role and we seek 
means to describe these also in the adelic setting. 
\\
We will realise non-holomorphic vector-valued Eisenstein series for the discrete Weil representation (cf.\ Definition~\ref{de:VVEisensteinseriesnonholomorphic}) in the present setting. 
To this end, let $(V,\qq)$ be the rational quadratic space associated to an even $\Z$-lattice $(L,\qq)$ 
and select the quadratic character $\chi = \chi_{V} = \Hsym_{\A}(\argdot,\disc(V))$. 
Recall that below~\eqref{eq:intertwiningoperator_finitepart}, 
for $\mu \in L'/L$ the finite adelic Schwartz form 
$\varphi_{\mu} \coloneq \mathds{1}_{\mu + \hat{L}} \coloneq \otimes_{p < \infty} \mathds{1}_{\mu + L_p}$ 
has been associated to a standard section $\Phi_{\mu}$ via the intertwining operator.  

\begin{de}\label{de:Eisensteinseries_VV_adelically}
	With the notation as above, we define the following \emph{Eisenstein series}  
	\begin{equation}
		\Eisadel_{\hat{L},l}(\tau,s) \coloneq \sum_{\mu \in L'/L} \Eisadel\left( \tau,s ; \Phi_{\infty}^{l} \otimes \Phi_{\mu}\right) \varphi_{\mu}. 
	\end{equation}
\end{de}

The shape of this series not only resembles that of the classical non-holomorphic vector-valued 
Eisenstein series presented in Definition~\ref{de:VVEisensteinseriesnonholomorphic} 
but may be directly translated to this setting \cite{BruinierYang2009}.  
For that purpose, an identification has to be established. 
Let $\Sspace_{L'/L}$ denote the subspace of $\Sspace(V \otimes \A_f)$ generated by 
the characteristic function \(\mathds{1}_{\mu + \hat{L}}\) for all elements $\mu \in L'/L$. 
Then this space is isometrically isomorphic to $\C[L'/L]$. 
It is also stable under a subrepresentation of the Weil representation of the group \(\Mp_{2}(\Z)\) 
whose dual is equivalent to the discrete Weil representation \(\dwrep_{L}\).
Additional details are provided in \cite{Bruinier2008}.

\begin{prop}\label{prop:Eisenstein_adelic_to_VVMF}
	With the notation as above we find the following identity 
	\begin{equation}\label{eq:prop:Eisenstein_adelic_to_VVMF}
		\Eisadel_{\hat{L},l}(\tau,s) = E_{L,0,l}\left(\tau,\tfrac{s+1-l}{2}\right).
	\end{equation}
	For this equality, the natural identification $\Sspace_{L'/L} \simeq \C[L'/L]$ is used. 
\end{prop}

\todo{This is lemma 3.9 in my thesis, but I should refer to Bruiniers and Yangs paper.}

\begin{proof}
	We follow \cite[Sec.~2]{BruinierYang2009} for part of the computation. 
	We begin our calculation by inserting the explicit term for the section $\Phi_{\infty}^{l}$ that has been calculated in Example~\ref{exa:standardsecweightlprincipalseriesinfty}. 
	\begin{align*}
			\Eisadel_{\hat{L},l}(\tau,s) 
		=\,	& \Im(\tau)^{-l/2}\sum_{\mu \in \Df} \Eisadel(g_{\tau} , s ; \Phi_{\infty}^l \otimes \Phi_{\mu}) \varphi_\mu \\
		=\,	&\frac{1}{2} \sum_{\mu \in \Df} \sum_{ \gamma \in \Gamma_{\infty} \backslash \Gamma(1) }  \frac{\Im(\gamma\tau)^{(s+1-l)/2}}{j(\gamma, \tau)^l} \Phi_{\mu}(\gamma) \varphi_\mu. 
	\end{align*}
	Now, recall that $\Phi_{\mu}$ was induced by $\lambda_{f} (\varphi_{\mu}) \in I_{f}(s_{0},\chi)$ 
	via the finite intertwining operator introduced in Example~\ref{exa:Principalseriesec_from_intertwining}. 	We then find via the correspondence of $\Sspace_{L'/L}$ and $\C[L'/L]$ that for $\gamma \in \Mp_{2}(\Z)$ 
	the following identity is true. 
	\begin{align*}
		\lambda_f(\varphi_\mu)(\gamma) 	
		&=	\langle \omega_f(\gamma) \varphi_{\mu} , \varphi_{0}  \rangle \\
		&=	\langle \overline{\omega_f(\gamma^{-1})} \varphi_{0} , \varphi_{\mu}  \rangle \\
		&=	\langle \rho_{L}(\gamma^{-1}) \varphi_{0} , \varphi_{\mu}  \rangle \\
		&=	\langle \rho_{L}(\gamma^{-1}) \efr_{0} , \efr_{\mu}  \rangle 
		=	\langle \efr_{0} , \rho_{L}(\gamma) \efr_{\mu}  \rangle. 
	\end{align*}
	Inserting this relation in the above expression for the Eisenstein series and using the isomorphism $\Sspace_{L'/L} \simeq \C[L'/L]$ we obtain 
	\begin{align*}
			\Eisadel_{\hat{L},l}(\tau,s) 
		\simeq\;	& \frac{1}{2} \sum_{ \gamma \in \overline{\Gamma_{\infty}} \backslash \Mp_{2}(\Z)}	 \frac{\Im(\gamma\tau)^{(s+1-l)/2}}{j(\gamma, \tau)^l}\sum_{\mu \in \Df}  \langle \rho_{L}(\gamma^{-1}) \efr_{0} , \efr_{\mu}  \rangle \efr_\mu \\
		=\;	& \frac{1}{2} \sum_{ \gamma \in \overline{\Gamma_{\infty}} \backslash \Mp_{2}(\Z)}	 \left[\Im(\tau)^{(s+1-l)/2}\right] \Big\vert_{l} \cdot \rho_{L}(\gamma^{-1}) \efr_{0} \\
				=\;	& E_{L,0,l}(\tau,\tfrac{s+1-l}{2}).
	\end{align*}
	This completes the proof. 
\end{proof}
\todo{Maybe insert a factor of $2$, depending on the convention for $\Gamma_{\infty}$ and $\overline{\Gamma_{\infty}}$ you choose.}

It is possible to compute the Fourier expansion of $\Eisadel(\tau,s; \Phi)$ for factorisable sections in terms of local Whittaker functions as explained in \cite[Thm.~2.4]{KudlaYangSL2}. 
However, the computations in concrete cases involve solving Gauss sums for the appearing local integrals and 
a sizable amount of case distinctions.

\subsection{Theta series}
In the classical setting for a positive definite quadratic space $(V,\qq)$ of even dimension \(m\), 
select the standard Gaussian \(\varphi_0\colon V_{\R} \to \C, \quad v \mapsto \exp(- 2 \pi \qq(v))\). 
Then for $\tau = u + iv \in \Ha$, $g_{\tau} = N(u) M(\sqrt{v})$ as in \eqref{eq:Iwasawa_decomposition} and $\omega_{\infty}$ the Weil representation we find 
\begin{equation}\label{eq:wrep_action_Schwartzform_infty}
	\omega_{\infty}(g_{\tau}) \varphi_0 (x) = e^{2 \pi i \qq(x) u} v^{\dimV/4} \varphi_0(\sqrt{v}x) = v^{\dimV/4} \exp( 2 \pi i \qq(x) \tau). 
\end{equation}
With this new perspective, we may define theta functions associated to Schwartz forms.

\begin{de}\label{de:thetafun_Schwartzform_coset}
	Let $\varphi \in \Sspace(V_{\R})$ and $\lambda \in \Df = L'/L$. 
	Then
	\begin{equation*}
		\theta_{L,\lambda}(\tau;\varphi) = v^{-\dimV/4} \cdot \sum_{x \in L + \lambda}  [\omega_{\infty}(g_{\tau}) \varphi](x)
	\end{equation*}
	is called \emph{theta function}\index{theta!function} associated to $\varphi$ and the coset $\lambda + L$.
\end{de}

We note that we may let $\OG(V_{\R})$ act on $\varphi$ via the left regular representation, i.e.\ 
\[
	h\colon \varphi \mapsto \varphi(h^{-1} \argdot) \in \Sspace(V_{\R}).
\]
In that sense, $\theta_{L,\lambda}$ may be understood as a function on $\Ha \times \OG(V_{\R})$.  
Combining the above theta functions for different cosets of the dual lattice yields the 
following theta function.

\todo{define on group and insert \(g_\tau\)}

\begin{de}\label{de:Thetafun_Schwartzform}
	Let $\varphi \in \Sspace(V_{\R})$. 
	Then
	\begin{equation*}
		\Theta_{L}(\tau,h;\varphi) = v^{-\dimV/4} \cdot \sum_{x \in L'}  [\omega_{\infty}(g_{\tau}) \varphi](h^{-1}x) \efr_{x + L} 
	\end{equation*}
	is called \emph{theta function}\index{theta!function} associated to $\varphi$.
	In case the Schwartz form $\varphi$ is invariant under the action of a maximal compact $K \subseteq \OG(V_{\R})$, the dependence on $h$ factors through $\Grassm$ and we write 
	$\Theta_{L}(\tau,z;\varphi)$ for $z \in \Grassm$. 
	Moreover, for fixed $z$, we omit this variable. 
\end{de}

A prominent example is the classical case based on \eqref{eq:wrep_action_Schwartzform_infty} 
given that $(V_{\R}, \qq)$ is positive definite. 
In the indefinite case, the naive Gaussian $V_{\R} \ni x \mapsto \exp(-2 \pi  \qq(x))$ is not rapidly decreasing and hence not a Schwartz function. 
For an element $z \in \Grassm$ we may decompose \(v = v_{z^{\perp}} + v_{z} \in z \oplus z^{\perp} = V\) and define the \emph{standard majorant} of \(\qq\) by   
\begin{equation}\label{eq:standardmajorant}
	\qq_{z}^{+}(v) = \qq(v_{z^{\perp}}) - \qq(v_{z})
\end{equation}
which is positive definite. 
This enables us to associate theta functions to indefinite quadratic spaces by introducing a parameter $z \in \Grassm$. 
In fact, with $\qq_{z}^{+}\colon V_{\R} \to \R$ as above, $\varphi(x) = e(i\qq_{z}^{+}(x))$, and $\tau = u + iv \in \Ha$  
we find 
\[
	\Theta_{L}(\tau,z;\varphi) = \sum_{x \in L'} e\left(u \qq(x) + iv \qq_{z}^{+}(x) \right) \efr_{x + L}. 
\]

Before finishing the current subsection we note that theta functions do factor, provided the associated Schwartz forms do.

\begin{rem}\label{rem:Thetafun_decomposes_tp_if_Schwartzform_does}
	If $L = L_1 \oplus L_2$ and $\varphi_{\infty} = \varphi_{\infty}' \otimes \varphi_{\infty}'' \in \Sspace(V_{1,\infty}) \otimes \Sspace(V_{2,\infty})$ split,
	then, provided the theta function $\Theta(\tau ; \varphi_{\infty})$ is absolutely convergent, 
	one has the following decomposition: 
	\begin{equation}
		\Theta_{L}(\tau; \varphi_{\infty}) = \Theta_{L_1}(\tau; \varphi_{\infty}') \otimes \Theta_{L_2} (\tau; \varphi_{\infty}''). 
	\end{equation}
\end{rem}

\begin{proof}
	A straight forward calculation yields the desired result 
	\begin{align*}
		\Theta_{L}(\tau ; \varphi_{\infty}) \;	
		=\;	&v^{-\frac{m}{4}} \cdot \hspace{-30pt} \sum_{(\lambda_1,\lambda_2) \in L_1'/L_1 \oplus L_2' / L_2} \sum_{l_1 \in L_1 + \lambda_1 , l_2 \in L_2 + \lambda_2} \left(	\omega_\infty (g_\tau) \varphi_{\infty} \left((l_1,l_2)\right)	\right) \e_{\lambda_1+\lambda_2}	\\
		=\;	&\sum_{\substack{\lambda_1 \in L_1'/L_1 \\ \lambda_2 \in L_2' / L_2}} \hspace{-10pt} v^{-\frac 1 4 } \hspace{-10pt} \sum_{l_1 \in L_1 + \lambda_1 } \hspace{-10pt} \left(\omega_\infty (g_\tau) \varphi_{\infty}'\right) (l_1)  \cdot v^{-\frac{m-1}{4}} \hspace{-10pt} \sum_{l_2 \in L_2 + \lambda_2} \hspace{-10pt}  \left(	\omega_\infty (g_\tau) \varphi_{\infty}'' \right) (l_2)	  \e_{\lambda_1+\lambda_2}	\\
		=\;	&\sum_{\substack{\lambda_1 \in L_1'/L_1 \\ \lambda_2 \in L_2' / L_2}} \theta_{L_1,\lambda_1} (\tau; \varphi_{\infty}') \cdot \theta_{L_2,\lambda_2} (\tau ; \varphi_{\infty}'')  \e_{\lambda_1+\lambda_2}	\\
		=\;	&\Theta_{L_1}(\tau; \varphi_{\infty}') \otimes \Theta_{L_2} (\tau; \varphi_{\infty}'').
	\end{align*} 
\end{proof}

		\section{Orthogonal setting}\label{sec:OrthogonalSetting}

Let $V$ be the rational quadratic space enveloping $(L,\qq)$ of signature $(\possig,2)$ and $V_{\R}$ its real points. 
The associated symmetric domain is the \emph{Grassmannian} 
\begin{equation}
	\Grassm(V_{\R}) \coloneq \set{ z \subseteq V_{\R} \mid \dim(z) = 2 \textnormal{ and } \qq\vert_{z} < 0 },
\end{equation}
where we understand the subspace to be oriented. 
It has two connected components and one of them will be denoted by \(\Grassm^+\). 
Note that by Witt's extension theorem, 
the group $\OG(V_{\R})$ acts transitively on $\Grassm$ and that there is a subgroup \(\OG^{+}(V_{\R})\) of index \(2\), 
acting transitively on \(\Grassm^{+}\). 
Further, we denote the kernel of the natural map 
\[
	\OG(L) \cap \OG(V_\R)^+ \to \Aut(L'/L)
\]
by \(\disckernel(L)\) and the associated modular variety by \(Y_L \coloneq \disckernel(L)\backslash \Grassm^{+}\). 

\nocite{Bruinier2008}

A convenient alternative realisation is the \emph{projective model}. 
To be more explicit, consider $V_{\C} = \C \otimes_{\Q} V$, which induces a projective variety 
	\[
		P(V_{\C}) = (V_{\C} \setminus \set{0}) / \C^{\times}.
	\]
Equipping this space with the final topology of the natural projection results in a compact space that  
naturally carries the structure of a complex analytic manifold. 
We consider the open subset of the zero quadric 
	\[
		\mathcal{K} \coloneq \set{ [Z] \in P(V_{\C}) \mid \bilf(Z,Z) = 0 , \bilf(Z,\overline{Z}) < 0}
	\]	
which carries a complex structure. 
The assignment $\mathcal{K} \ni [Z] \mapsto \R \Re(Z) + \R \Im(Z) \in \Grassm^{+}$ defines a continuous and open two-to-one map, 
transferring the complex structure. 

We choose a nonzero isotropic vector $e_{1} \in V$ and some $e_{2} \in V$ 
with $\bilf(e_{1},e_{2}) = 1$. 
The space $W = e_{1}^{\perp} \cap e_{2}^{\perp}$ is of type $(\possig-1,1)$, 
so that 
\begin{align}\label{eq:tubedomainmodel_Vdecomp}
	V = W \oplus \Q e_{1} \oplus \Q e_{2}.
\end{align}
Then we define 
\begin{equation}\label{eq:Tubedomain_full}
	\Ho_{\possig}^{\pm} \coloneq \{ z = x + i y \in W_{\C} \mid \qq(y) < 0 \}.
\end{equation}
We write $Z=(z,a,b) \in V_{\C}$ for an element in the associated complex space given by $z + a e_{1} + b e_{2}$ 
with $ z \in W_{\C}$. Then the map 
\[
	\Ho_{\possig}^{\pm} \ni z \mapsto \psi(z) = [(z,-\qq(z)-\qq(e_{2}),1)] \in \mathcal{K}
\]
defines a biholomorphic association. 
Set $\mathcal{H}_{\possig}^{+}$ and $\mathcal{K}^{+}$ to be the respective connected components in bijection with \(\Grassm^{+}\). 
The action of \(\OG(V_\R)\) on \(P(V_{\C})\) carries over to \(\mathcal{H}_{\possig}^{+}\) and induces an action by fractional linear transformations. 
This yields for \(\gamma \in O(V_\R)\) a \emph{factor of automorphy} \(j(\gamma,z)\) such that \(\gamma \psi(z) = j(\gamma,z) \psi(\gamma z)\). 
Then for $r \in \Q^{\times}$ there is a map $w_r$ on $\OG^{+}(V_{\R}) \times \OG^{+}(V_{\R})$ to the roots of unity of order bounded by $\den(r)$ such that 
\[
	 j(\gamma_1 \gamma_2, Z)^{r} = w_r(\gamma_1,\gamma_2) \cdot j(\gamma_1,\gamma_2z)^{r} j(\gamma_2,z)^{r}.
\]

Similar to the elliptic setting there is a notion of modular forms in the orthogonal case. 
We refer to \cite[3.3 p.~84]{Bruinier2002} for a more detailed description. 

\begin{de}
	Let $\Gamma \leq \OG^{+}(V_{\R})$ and $r \in \Q$. 
	A multiplier system\index{multiplier system} of weight $r$ for $\Gamma$ is a map
	\[
		\chi\colon \Gamma \to \Tc, \qquad \textnormal{such that} \qquad \chi(\gamma_1 \gamma_2) = w_r(\gamma_1,\gamma_2)^{-1} \cdot \chi(\gamma_1)\chi(\gamma_2). 
	\]
\end{de}

With the notion of a multiplier system, modular forms on the tube domain $\Ho_{\possig}^{+}$ may be defined. 
Let $V$ be isotropic. 
An orthogonal \emph{meromorphic modular form}\index{meromorphic!modular form}\index{modular form!meromorphic}\index{modular form!orthogonal} with respect to the group 
$\Gamma \leq \disckernel[L]$ of finite index of weight $r \in \Q$ and multiplier system $\chi$ 
is a meromorphic function $F$ on $\Ho_{\possig}^{+}$ with transformation property 
\begin{equation*}
	F(\gamma z) = \chi(\gamma) j(\gamma ,z)^{r} F(z).
\end{equation*}

When defining the Petersson slash operator in analogy to the elliptic case, 
Eisenstein series may be constructed in a similar fashion as an example of modular forms.
We will not require these in the following, but the curious reader may consult \cite[7 p.~2867]{Kiefer2022}. 
Further, non-trivial examples of modular forms may be constructed via the Borcherds lift, 
relating weakly holomorphic vector-valued elliptic modular forms to meromorphic orthogonal modular forms (see\ Section~\ref{sec:ThetaLifts}).

\subsection{Special divisors}\label{ssec:Special_divisor}

In algebraic geometry divisors are a common tool for investigating the structure of varieties. 
A feature of the orthogonal setting is that the associated Shimura varieties have an abundance of divisors 
being induced by sub-varieties of the same type in all codimensions. 
These give rise to so called \emph{special cycles} which 
will play a key role in proving the injectivity of the Kudla--Millson lift in Section~\ref{sec:KMlift_Injectivity}. 
In the following, we describe the construction of these special cycles of codimension $1$.

Recall that $(V,\qq)$ has signature $(\possig,2)$ 
and select some $v \in V$ with $\qq(v) > 0$. 
Then the subspace $V_{v} \coloneq v^{\perp}$  \todo{better start with subspace ... element of $\Grassm$?} 
has signature $(\possig-1,2)$ and 
\begin{align}\label{eq:de_canonical_divisor}
	\candiv_{v} \coloneq \set{[z] \in \mathcal{K}^{+} \mid (z,v) = 0}
\end{align}
defines an analytic divisor on $\mathcal{K}^{+}$ being exactly the hermitian symmetric domain corresponding to~$(V_{v} , \qq)$. 
In other words, we find  
\begin{equation}\label{eq:subgrassmannian}
	\candiv_{v} \simeq \Grassm_{v} \coloneq \{ U \subset V_{v}(\R) \mid \dim(U) = 2 \textnormal{ and } \qq\vert_{U} < 0 \}.
\end{equation}
Here, $\Grassm_{v}$ can be realised as the Grassmannian attached to the stabiliser $\OG(V_{\R})_{v}$ of $v$ in $\OG(V_{\R})$. 
We call $\candiv_{v}$ a \emph{special divisor}\index{divisor!special}. 
Its corresponding description in the tube domain model looks as follows. 
Write $v = v_{W} + a e_{2} + b e_{1}$ as above (cf.\ \eqref{eq:tubedomainmodel_Vdecomp}) 
so that 
\[
	\candiv_{v} \simeq \set{z \in \mathcal{H}^{+} \mid a\qq(z) - \bilf(z,v_{W}) - a\qq(e_{2}) - b = 0} \subset \mathcal{H}^{+}.
\]
This description gives rise to the term \emph{rational quadratic divisor}\index{divisor!rational quadratic}. 
Further, the following combination of these divisors appears naturally in numerous applications (cf.\ Theorem~\ref{thm:Borcherds_theorem}). 

\begin{de}\label{de:Heegner_divisor}
	Select $\lambda \in L'/L$ as well as a rational number $n > 0$. 
	Then 
	\begin{equation}\label{eq:de:Heegner_divisor}
	\Heegnerdiv(\lambda , n) \coloneq \sum_{\substack{ v \in \lambda + L \\ \qq(v) = n}} \candiv_{v}
	\end{equation}
	defines an analytic divisor on $\mathcal{K}^{+}$ which is called \emph{Heegner divisor}\index{Heegner!divisor} of discriminant $(\lambda,n)$. 
\end{de}

Recall that to an arithmetic subgroup $\Gamma \leq \OG(L)$ there is 
an (orthogonal) Shimura variety 
$Y_{\Gamma} = \Gamma \backslash \Grassm$. 
Now, if $\Gamma$ lies in the discriminant kernel, 
i.e.\ acts trivially on the discriminant group $L'/L$, 
then the Heegner divisor 
descends to an algebraic divisor on $Y_{\Gamma}$, which is also denoted by~$Z(\lambda,n)$. This may be verified by Chow's lemma. 
\todo{take a look at Ku2 from JHB123 to verify defined over number field}

Alternatively, the Heegner divisors might be symmetrised over $L'/L$ to yield an algebraic divisor on $Y_{\Gamma}$, namely 
\begin{equation}
	\Heegnerdiv(m) \coloneq \frac 1 2 \sum_{\lambda \in L'/L} Z(\lambda,m).
\end{equation}

We also require the following more primitive divisor which is directly induced by \eqref{eq:de_canonical_divisor}. 
Let $\Gamma_{v} \coloneq \Gamma \cap \OG(V_{\R})_{v}$ be the intersection with the stabiliser of $v$. 
Then we find that 
\begin{equation}\label{eq:de:candiv}
	\candiv(v) \coloneq \Gamma_{v} \backslash \Grassm_{v} \to Y_{\Gamma}
\end{equation}
defines a (in general relative) cycle. These cycles play a fundamental role in the proof of the injectivity of the Kudla--Millson lift 
in Section~\ref{ssec:Cycle_Integrals_KM_Lifts}. 

The existence of such a family of algebraic divisors is distinctive for orthogonal and unitary groups.

		\section{Siegel-Weil formula}\label{sec:SiegelWeil}

Theta series and Eisenstein series form two fundamental families of automorphic forms. 
In many cases
(cf.\ \cite[Thm.~6.7]{Mueller2024}) 
the space of cusp forms is spanned by theta series, which is complemented by Eisenstein series. 
A central tool relating theta series to Eisenstein series is the \emph{Siegel--Weil} formula.
In \ref{ssec:geometric_interpretation} we follow Kudla \cite[Sec.~4.3]{Kudla2003} 
in applying it to compute certain geometric integrals of the Kudla--Millson Schwartz form, which will be essential for the injectivity arguments in Section~\ref{sec:KMlift_Injectivity}.
Recall that $L$ is an even lattice of signature $(\possig,2)$ and the dimension of the rational quadratic space $V = L \otimes_{\Z} \Q$ 
is denoted by $\dimV$.

\subsection{Siegel--Weil formula}\label{ssec:SWF}

The subject of this subsection is the so called \emph{Siegel--Weil} formula which establishes a relation between Eisenstein series as in \eqref{eq:EisensteinSeries} 
evaluated at a special point $s_{0} \in \C$ and theta integrals as described in \eqref{eq:ThetaIntegral}. 
It is an adelic extension of the classical Siegel formula and has first been proven by André Weil \cite{WeilSiegel}.

For a Schwartz--Bruhat function $\varphi \in \Sspace(V_{\A})$ define the adelic theta distribution 
\[
	\theta(g',h;\varphi) = \sum_{x \in V_{\Q}} (\omega(g')\varphi) (h^{-1}x)
\]
as a function in $g' \in \Mp_{2}(\A)$ and $h \in \OG(V_{\A})$. 
Averaging over the argument $h$ delivers the following standard theta integral.

\begin{de}\label{de:Theta_integrals}
	For a Schwartz--Bruhat function $\varphi \in \Sspace(V_\A)$ 	define the following theta integral.  
	\begin{equation}\label{eq:ThetaIntegral}
		I(g';\varphi) \coloneq \int_{\OG(V_\Q) \setminus \OG(V_\A)} \theta(g',h;\varphi) \di h,
	\end{equation}
	where $\textrm{d} h$ denotes the standard invariant normalised measure on the quotient such that $\vol(\OG(V_\Q) \setminus \OG(V_\A)) = 1$.
\end{de}

By Weil's convergence criterion\index{Weil!convergence criterion}, the integral $I(g;\varphi)$ above is absolutely convergent in case the rational space $V$ is anisotropic or $\dim(V) = \dimV > 2 + r$, 
where $r$ denotes the Witt rank of $V$. 
In case this criterion is not matched, the integral may still be regularised 
as in \cite[p.~41]{KudlaRallis1994} and numerous results may be transferred to this setting. 
With these preparations we require only one additional piece of notation to state the Siegel--Weil formula. 
Let $\varphi \in \Sspace(V_{\A})$ be a Schwartz--Bruhat form and $\Phi \in I(s,\chi)$ its associated 
standard section through the principal series representation arising from the intertwining operator 
as described in~\eqref{eq:intertwiningoperator_finitepart} (also compare Lemma~\ref{lem:Intertwiningintoprincseries}). 

\begin{de}\label{de:Eisenstein_adelically_depence_Schwartzform}
	With the notation as above, write 
	\[
		\Eisadel(g',s;\varphi) \coloneq \Eisadel(g',s;\Phi), \qquad \Eisadel(\tau,s; \varphi) \coloneq \Eisadel(\tau,s;\Phi). 
	\]
\end{de}

\todoexp{IntBorchForms p.325 (4.9)-- Find ref and explain properly (above in Eisenstein section)}

\begin{thm}[Siegel--Weil formula, {\cite[Thm.~4.1]{Kudla2003}}]\label{thm:SiegelWeil}
	Assume $\dimV > 2 + r$ or that $V$ is anisotropic. 
	Then the theta integral $I(g';\varphi)$ 			for $\varphi \in \mathcal{S}(V_{\A})$ is absolutely convergent. 
		Further, $\Eisadel(g',s;\varphi)$ is holomorphic at $s_{0}= \frac{\dimV}{2} - 1$ and the following identity holds: 
		\[
			\Eisadel(g',s_{0};\varphi) = \kappa \cdot I(g';\varphi),
		\]
	where $\kappa = 2$, if $\dimV \leq 2$, and $\kappa = 1$ otherwise. 	
\end{thm}

In particular, the Eisenstein series associated to a Schwartz--Bruhat form bears information 
about representation numbers of a quadratic lattice in its Fourier expansion at the critical point $s_{0} = \frac{\dimV}{2} - 1$.  
This is remarkable since the Eisenstein series is essentially built from local data, 
while the theta integral involves global arithmetic.
Hence, theta integrals associated to different quadratic spaces 
may be related to each other, via the Eisenstein series, once the local data agrees.

\begin{de}
	Let $V_{p},V_{p}'$ be two quadratic $\Q_{p}$ spaces of dimension $\dimV$. 
	Suppose these spaces induce the same quadratic character $\chi_{V_p} = \chi_{V_p'}$ via the Hilbert symbol. 
	Two functions $\varphi_{p} \in \Sspace(V_{p})$ and $\varphi_{p}' \in \Sspace(V_{p}')$ are called \emph{matching}\index{function!matching}, if they give rise to the same section via the intertwining operator $\lambda_{p}(\varphi_{p}) = \lambda_{p}'(\varphi_{p}')$. 
\end{de}

In fact, in case of $\dimV > 4$ and $p$ non-archimedean, every $\varphi_{p} \in \Sspace(V_{p})$ has a 
matching function~$\varphi_{p}' \in \Sspace(V_{p}')$, since the 
local principal series $I_{p}(s_0,\chi_{p})$ is irreducible. 

In order to compare theta integrals as suggested above, it is hence sufficient to find a non-trivial match 
at the archimedean place. One solution to that problem is given by the Kudla--Millson Schwartz function 
which is introduced in the following section.

\subsection{Kudla--Millson Schwartz form}\label{ssec:KMform-properties}

Let $(V,\qq)$ have signature $(\possig,2)$ and dimension~$\dimV$. 
Recall that for a complex manifold $D$, the space of smooth forms of Hodge type $(1,1)$ 
is denoted by~$\Dforms^{(1,1)}(D)$.
Kudla and Millson have constructed a form $\KMfun \in \Sspace(V_\R) \otimes \Dforms^{(1,1)}(\Grassm)$, 
where~$\Grassm$ denotes the oriented Grassmannian associated to the base lattice $L$, with the following properties.

\begin{enumerate}[label = \roman*)]
	\item 
	For all $h \in \OG(V_\R)$ 
	\begin{equation}\label{eq:KMform_Oinvariance}
		h^\ast \KMfun(h^{-1}x) = \KMfun(x).
	\end{equation}
	
	\item 
	The form $\KMfun$ has weight $\dimV/2$ for $K_\infty'$ for the Weil representation $\omega_{\infty}$, 
	meaning 
	\begin{equation}
		\omega_\infty(k')\KMfun = \nu_{\dimV/2}(k') \KMfun, 
	\end{equation}
	where the character $\nu_{l}$ is presented in Remark~\ref{rem:MSO_characters}.
	
	\item 
	The form $\KMfun$ is closed, i.e.\ 
	\begin{equation}
		\di \KMfun = 0
	\end{equation}
	for the exterior differential $\di$ on $\Grassm$. 
\end{enumerate}

This form \emph{matches} the Gaussian on a positive definite space of the same dimension in the sense that it gives rise to the same section through the principal series representation. 
To specify the meaning of this statement, we define the following form.

\begin{de}\label{de:Omega_form}
	Define the following closed $\OG(V_{\infty})$ invariant $(1,1)$ form on $\Grassm$:
	\begin{equation}
		\Kaehlerform \coloneqq \KMfun(0).
	\end{equation}
\end{de}

An explicit description in local coordinates 
on the tube domain model of that form is presented in \cite[Prop.~4.11]{Kudla2003}. 
It is given by 
\begin{equation}\label{eq:Omega_in_coordinates}
	\Kaehlerform = - \frac{1}{2 \pi i} \left[ - \bilf(y,y)^{-2} \bilf(y,\di z) \wedge \bilf(y,\di \overline{z}) + \bilf(y,y)^{-1} \tfrac{1}{2} \bilf(\di z , \di \overline{z}) \right]
\end{equation}
in coordinates on the tube domain model $\Ho_{\possig}^{\pm}$ (cf.\ \eqref{eq:Tubedomain_full}).
Its negative $-\Kaehlerform$ is an invariant \emph{Kähler form}\index{Kähler!form} on $\Grassm$ and descends to a Kähler form on the variety $X_{K}$ given in \eqref{eq:Adelic_Spin_variety}.

Note that we may write 
\[
\KMfun(x) \wedge \Kaehlerform^{\dimV-1} = \KMfunt(x)\Kaehlerform^{\dimV}
\]
for a function $\KMfunt \in \Sspace(V_{\R}) \otimes A^{(0,0)}(\Grassm)$. 
The $\OG(V_{\infty})$ invariance of $\Kaehlerform$ then implies that 
\[
\KMfunt(hx,hz) = \KMfunt(x,z)
\]
for any $h \in \OG(V_{\infty})$. 
Further, setting $k \coloneq \tfrac{\dimV}{2}$ we find 
\begin{align}\label{eq:KMfun_intertwines_to_Phik}
	\lambda_{\infty}(\KMfun(\argdot,z)) = \Phi^{k}(s_0) \Kaehlerform, \quad \lambda_{\infty}(\KMfunt(\argdot,z)) = \Phi^{k}(s_{0})
\end{align}
for any $z \in \Grassm$, where $\lambda_{\infty}$ denotes the intertwining operator described in 
Lemma~\ref{lem:Intertwiningintoprincseries} and $\Phi_{\infty}^{k}$ denotes the standard section 
in $I_{\infty}(s,\chi_{V})$ at the infinite place presented 
in Example~\ref{exa:standardsecweightlprincipalseriesinfty}.  
This means, the function $\KMfunt$ will be associated by $\lambda_{\infty}$ to the same section as 
the standard Gaussian $\varphi_{0}$ on a space $V'$ of signature $(\dimV,0)$ if the local characters agree, i.e.\ they are \emph{matching}. 
If further, we have two matching functions at the finite places $\varphi_{f} \in \Sspace(V_f)$ and $\varphi_{f}' \in \Sspace(V_{f}')$, the Siegel--Weil formula (Theorem~\ref{thm:SiegelWeil}) yields 
\begin{equation}
	I(g'; \KMfunt \otimes \varphi_{f}) = E(g', s_0; \KMfunt \otimes \varphi_{f}) = E(g',s_0; \varphi_{0} \otimes \varphi_{f}') = I(g'; \varphi_{0} \otimes \varphi_{f}').
\end{equation}

Before turning towards an essential application of the Siegel--Weil formula within the scope of this paper, 
we remark that a more explicit representation of the Kudla--Millson Schwartz form 
is found in \cite{Metzler2026}.

\subsection{A geometric integral}\label{ssec:geometric_interpretation}

Denote the GSpin group 
associated to the lattice $(L,\qq)$ of signature $(\possig,2)$ by $H$. 
Moreover, denote the oriented Grassmannian of $L$ by $\Grassm$ and
choose a compact open subgroup $K_{f} \subset H(\A_{f})$.
Then we consider the variety 
\begin{equation}\label{eq:Adelic_Spin_variety}
	X_{K} \coloneq H(\Q) \backslash (\Grassm \times H(\A_{f})) / K_{f}
\end{equation}
and have that 
\[
	X_{K} = \bigsqcup_{j = 1} \Gamma_{j} \backslash \Grassm^{+},
\]
where $\Grassm^{+}$ is a connected component of $\Grassm$. 
We write $Y_{j} \coloneq \Gamma_{j} \backslash \Grassm^{+}$ and note that there are
$\lvert \hat{\Z}^{\times} : \Cliffnorm(K_{f})\rvert$ of these components. 
Here, $\Cliffnorm$ denotes the Clifford norm\index{Clifford!norm}. 
With this preparation we may state the following modified version of \cite[Prop.~4.17]{Kudla2003}.

\begin{prop}\label{prop:Geometric_Theta_integral}
	We find for a compact subgroup $K_f \leq H(\A_f)$ with $K_f \cap Z(\A) \simeq \hat{\Z}^\times$ under the isomorphism $Z(\A) \simeq \A^{\times}$ that 
	\[
		(-1)^{\possig} \tfrac{1}{4} \vol(K_{f}) \cdot \int_{X_K} \theta(g';\KMform \otimes \varphi_{f}) \wedge \Kaehlerform^{\possig - 1} = I(g'; \KMfunt \otimes \varphi_{f}). 
	\]
	Further, the integral on the right hand side over each component $Y_{j}$ of $X_{K}$ is identical. 
	Hence, for any component $Y_{j}$ the following identity holds 
	\[
		(-1)^{\possig} \tfrac{1}{4} \vol(K_{f}) \cdot \lvert \hat{\Z}^{\times} : \Cliffnorm(K_{f})\rvert \cdot \int_{Y_j} \theta(g';\KMform \otimes \varphi_{f}) \wedge \Kaehlerform^{\possig - 1} = I(g'; \KMfunt \otimes \varphi_{f}). 
	\]
\end{prop}

We would like to apply this result in the more classical setting of varieties defined at the infinite place over 
the orthogonal group~$\OG$ as in \cite{Borcherds1998} and recall that
\[
	1 \to \Q^{\times} \to H(\Q) \to \SO(V) \to 1,
\]
as well as 
\[
	H(\hat{\Z}) \to \disckernelso[\hat{L}] \to 1,
\]
where $\disckernelso[\hat{L}]$ denotes the discriminant kernel in $\SO(\hat{L})$. 
For $\Gamma \leq \SO(L)$, we may consider the closure $\overline{\Gamma} < \SO(V_{\A_{f}})$ 
and further select a cover $K_{f} < H(\A_{f})$, containing $\hat{\Z}^{\times}$. 
In that case, we have that the Shimura variety $X_{K}$ of the pair $(H,K_{f})$ is 
the same as the variety associated to $(\SO,\overline{\Gamma})$ (cf.\ \cite[p.~1669]{Bruinier2021}): 
\[
	\SO(V) \backslash \Grassm \times \SO(V_{\A_{f}}) / \overline{\Gamma} \simeq H(\Q) \backslash \Grassm \times H(\A_{f}) / K_{f}.
\]
In the special case of $K_{f} = \GSpin_{\hat{L}}$, 
we find that for the choice of $h_{j}$ being equal to the neutral element that 
$Y_{j} \simeq Y'(L) = \disckernelso[L] \backslash \Grassm^{+}$, where $\disckernelso[L]$ denotes the discriminant kernel in $\SO$. 
In fact, 
if we specialise further, to the case of $L$ splitting a hyperbolic plane, 
we find that the index $\lvert \hat{\Z}^{\times} : \Cliffnorm(K_{f})\rvert$ equals $1$, 
so that there is only one connected component identifying with $Y'(L) = \disckernelso[L] \backslash \Grassm$. 
All in all, we obtain the following corollary. 

\begin{cor}\label{cor:OCycleint_as_standardtheta}
	Let $(L,\qq)$ be an even quadratic lattice of signature $(\possig,2)$ 
	and select $K_{f} = \GSpin_{\hat{L}}$. 
	Further, write $Y_{L} = \disckernel[L] \backslash \Grassm$ 
	for the variety associated to the discriminant kernel $\disckernel[L]$ in $\OG(L)$ 
	and select an element $\varphi_{f} \in \Sspace(V_{\A_{f}})$ that is invariant under $\disckernel[L]$. 
	Then 
	\begin{equation}\label{eq:cor:OCycleint_as_standardtheta}
		\int_{Y_{L}} \theta(g';\KMfun \otimes \varphi_{f}) \wedge \Kaehlerform^{\possig - 1} 
		= 
		(-1)^{\possig} \frac{4 c_{K}}{\lvert \hat{\Z}^{\times} : \Cliffnorm(K_{f})\rvert} \cdot 
		I(g',\KMfunt \otimes \varphi_{f}),
	\end{equation}
	where $c_{K} = 1$ or $1/2$, depending on whether $\disckernel[L] \cap \SO(V) \backslash \Grassm \to \disckernel[L] \backslash \Grassm$ is bijective or two-to-one. 
\end{cor}

\begin{proof}
	If considered in the special orthogonal setting, meaning integrating over 
	$\disckernelso[L]\backslash \Grassm$ instead of $\disckernel[L]\backslash \Grassm$ on the left hand side 
	of \eqref{eq:cor:OCycleint_as_standardtheta}, the statement follows immediately from 
	Proposition~\ref{prop:Geometric_Theta_integral}. 
	The only task remaining is hence the translation to the orthogonal setting. 
	By assumption, the integrand is invariant under pullback by elements from $\disckernel[L]$. 
	In fact, the form $\Kaehlerform$ is invariant under $\OG(V_{\infty})$ and $\KMfun$ fulfils \eqref{eq:KMform_Oinvariance}. 
	Further, $\varphi_{f}$ was chosen to be invariant under $\disckernel[L]$ as well, 
	so that we may translate the domain of integration by 
	$\gamma \in \disckernelso[L] \backslash \disckernel[L]$ and combine a fundamental domain for 
	$\disckernelso[L] \backslash \Grassm$ from translates of a domain for 
	$\disckernel[L]\backslash \Grassm$. 
	Consequently, the integral in question in \eqref{eq:cor:OCycleint_as_standardtheta} equals the integral over 
	$\disckernelso[L] \backslash \Grassm$ divided by the multiplicity of the cover $\disckernelso[L] \backslash \Grassm \to \disckernel[L]\backslash \Grassm$. 
	\\
	Finally, we note that
	\[
		[\disckernel[L] : \disckernelso[L]] \leq [\OG(V) : \SO(V)] = 2, 
	\]
	meaning that the considered cover represents at most a double-cover. 
\end{proof}

		\section{Theta lifts}\label{sec:ThetaLifts}

The two theta lifts central to our considerations are the Borcherds lift and the Kudla--Millson lift.
The former provides a multiplicative correspondence from vector-valued weakly holomorphic modular forms to orthogonal modular forms with explicit product expansions. 
The latter constitutes a linear lift from cusp forms to harmonic differential forms on orthogonal Shimura varieties.
These two lifts are closely related: in particular, injectivity of the Kudla--Millson lift is, under suitable hypotheses, equivalent to a converse theorem for Borcherds products (see Theorem~\ref{thm:KMliftinj_Borliftsurj}).
In this section we recall their definitions and review known injectivity results.
Let \(L\) denote an even lattice of signature \((\possig,\negsig)\).

\subsection{The Borcherds lift}

In his 1995 paper \cite{Borcherds1995}, Borcherds introduced a multiplicative lift from weakly holomorphic modular forms of weight $1/2$
to meromorphic modular forms, producing infinite product expansions. 
The weights and divisors of the resulting forms are determined explicitly in terms of the Fourier coefficients of the input.
Notable applications include infinite product expansions for the discriminant function, the \(j\)-invariant, and the Monster denominator formula 
used to prove the moonshine conjecture. 

In his renowned 1998 paper \cite{Borcherds1998} Borcherds significantly expanded on his initial approach, 
presenting a multiplicative lift that maps weakly holomorphic 
vector-valued modular forms of weight $1-\possig/2$ 
to meromorphic automorphic forms for certain arithmetic subgroups of $\OG_{\possig,2}(\R)$, 
i.e.\ to the orthogonal setting. 
In the case of $\possig = 1$, the association may be interpreted as a purely symplectic lift, 
due to the accidental isomorphism of $\SO_{1,2}(\R)$ and $\Sp_{1}(\R)$.
For a concise source on the construction  
see \cite{BruinierConverse}.

\begin{thm}[Borcherds]\label{thm:Borcherds_theorem}
	Let $f \in \MF_{L^-,1-\possig/2}^{!}$ be a weakly holomorphic modular form with Fourier coefficients $a(\lambda,n)$ as in \eqref{eq:VVMF_FourierExp}. Assume that $a(\lambda,n)$ is integral for $n \leq 0$. 
	Then there exists a meromorphic modular form $\Psi(\argdot,f)$ for the discriminant kernel $\disckernel[L]$ 
	with unitary multiplier system of finite order, such that the following conditions are satisfied. 
	\begin{enumerate}[label=\alph*)]
		\item 
		The weight of $\Psi(\argdot,f)$ is given by $a(0,0)/2$.
		
		\item 
		The divisor of the target form $\Psi( \argdot , f)$ is determined by the principal part of $f$: 
		\begin{equation}\label{eq:Psi_divisor}
			\frac{1}{2} \sum_{\mu \in L'/L}\sum_{n > 0} a(\mu,-n) Z(\mu,n),
		\end{equation}
		where 
		\[
		Z(\mu,n) = \sum_{\substack{v \in \mu + L \\ \qq(v) = n}} \{z \in \Grassm \mid z \perp v\}
		\]
		is the \emph{Heegner divisor}\index{Heegner!divisor} of discriminant $(\mu, n)$ (cf.\ Definition~\ref{de:Heegner_divisor} plus the subsequent comment). 		
		\item 
		The target form $\Psi(\argdot,f)$ has explicit infinite product expansions at cusps. 
	\end{enumerate} 
\end{thm}

For further details regarding the product expansions that are constructed for different cusps, 
or background material on the construction itself, 
we refer the reader to \cite[Thm.~13.3]{Borcherds1998} and conclude with the following definition.

\begin{de}
	For a form $f \in \MF_{L^-,1-\possig/2}^{!}$ as above, we refer to $\Psi(\argdot,f)$ as the \emph{Borcherds lift}\index{Borcherds!lift}\index{lift!Borcherds} of the form $f$ or, alternatively, as the \emph{Borcherds product}\index{Borcherds!product}\index{product!Borcherds} of $f$. 
\end{de}

\subsection{The Kudla--Millson lift}\label{ssec:KMlift}

Recall that the Kudla--Millson Schwartz form $\KMform$ has been introduced in \ref{ssec:KMform-properties}. 
It is noteworthy that $\Kaehlerform \coloneq \KMform(0, \argdot)$ is the Euler form 
of the symmetric space~$\Grassm$ 
and in the hermitian case of $\negsig = 2$, its negative defines a Kähler form. 
In the latter case, a description in coordinates on the associated tube domain model has been provided in~\eqref{eq:Omega_in_coordinates}. 

With the Kudla--Millson Schwartz form available, the associated theta kernel may be constructed as presented 
in Definition~\ref{de:Thetafun_Schwartzform}. 

\begin{de}\label{de:KMTheta}
	For $\tau = u + iv \in \Ha$ and $z \in \Grassm$ the \emph{Kudla--Millson theta function}\index{KudlaMillson@Kudla--Millson!theta function} of the lattice $L$ is given by 
	\begin{equation}
		\KMTheta(\tau,z) \coloneq \Theta_{L}(\tau,z;\KMform) = v^{-\dimV/4} \sum_{\lambda \in L'/L} \sum_{l \in \lambda + L} \left( \wrep_{\infty}(g_{\tau}) \KMfun\right) (l,z) \efr_{\lambda}. 
	\end{equation}
\end{de}

It should be noted that there is a more explicit construction provided in \cite[Sec.~2,4]{Bruinier2004a} based on which the
form is explicitly expressed in terms of Siegel theta functions in~\cite{Metzler2026}. 
As with the Borcherds lift, the theta kernel associated to the Kudla--Millson Schwartz form may be 
utilised as an integral kernel to construct a theta lift.

\begin{de}\label{de:KMlift}
	Let $L$ be an even lattice of signature $(\possig, \negsig)$ with rank $\dimV$ and
	set $k \coloneq \dimV/2$. 
	Then for a cusp form $f \in \CF_{L,k}(\Mp_{2}(\Z))$ and $\tau = u + iv \in \Ha$, 
	the association 
	\begin{equation}
		f \mapsto \KMlift(f) \coloneq \int_{\Mp_{2}(\Z) \backslash \Ha} v^{k} \langle f(\tau) , \Theta_{L}(\tau,z;\KMform) \rangle \, \frac{\di u \di v}{v^{2}}
	\end{equation}
	defines a map to the $\negsig$-forms $\Dforms^{\negsig}(\Grassm)$ 
	and is referred to as the \emph{Kudla--Millson lift}.\index{KudlaMillson@Kudla--Millson!lift}\index{lift!Kudla--Millson}
\end{de}

The target form $\KMlift(f)$ is in general closed and inherits the invariance under 
the discriminant kernel subgroup $\disckernel[L]$ of $L$ on the orthogonal side from the theta form~$\KMTheta$. 
Consequently, it descends to a form on the quotient space $Y_{L} = \disckernel[L] \backslash \Grassm$. 
In the hermitian case of~$\negsig = 2$, we obtain a map
\begin{equation}\label{eq:de:KMlift}
	\KMlift: \CF_{L,1+\possig/2} \to \HF^{(1,1)}(Y_{L}), 
\end{equation}
to the space of square-integrable harmonic differential forms of Hodge type $(1,1)$ 
on the orthogonal modular variety associated to the discriminant kernel $\disckernel[L]$. 

The applications of the Kudla--Millson theta correspondence range from the study of the cohomology 
of orthogonal and unitary Shimura varieties \cite{Kudla1986} and Arakelov theory of Shimura varieties \cite{Kudla2004a} over specific counting problems \cite{Engel2023},  
constructing mock modular forms and higher dimensional error functions \cite{Funke2017}, to proving a converse theorem for Borcherds products \cite{Bruinier2002} \cite{BruinierConverse}. 
It is the latter application that we have in mind which is related to the injectivity of the Kudla--Millson lift. 
The question of its injectivity already arose in \autocite{Kudla1990} 
and may be used to compute the rational Picard number of the underlying Shimura variety \autocite{Bergeron2016}, 
as well as to derive properties of cones generated by special cycles \autocite{Bruinier2019} \autocite{Zuffetti2022}. 
Since the Kudla--Millson lift has been constructed, there have been a number of advances in proving its injectivity. 
To the best of the author's knowledge, the first such result in the literature is the following.

\begin{thm}[{\cite[Thm.~5.12]{Bruinier2002}}]
	Let $(L,\qq)$ be an even lattice of signature $(\possig,2)$ such that it splits two hyperbolic planes $L \simeq K \oplus H \oplus H$. 
	Then the Kudla--Millson lift $\KMlift : \CF_{L,1 + \possig/2} \to \Hforms^{(1,1)}(Y_{L})$ is injective. 
\end{thm}

The proof is based upon the computation of the Fourier expansion of the target form, 
which reveals information about the Fourier coefficients of the initial form $f$ and leads to the conclusion 
that it must vanish if its image does. 
The second hyperbolic split is required in order to guarantee that indices of Fourier coefficients are 
representable by the lattice in a certain fashion similar to Lemma~\ref{lem:Hypsplit_Indices_primitivelyrep}. 

There has been a second advance by the same author in collaboration with Jens Funke in \cite{BruinierFunkesurj2010}. 
In this instance, the authors develop a new strategy of proof, 
by utilising the doubling method to compute the $L^2$--norm of 
the Kudla--Millson lift by means of the Rallis inner product formula 
in order to conclude injectivity in the unimodular case for general signature.  
In addition, the authors consider a twist of the Kudla--Millson lift, 
denoted by~$\KMFlift{l}$, which was introduced by 
Funke and Millson \cite{Funke2012} with a parameter $l \in \N_0$ 
that allows for inputting forms of higher weight. 
The twisted lift $\KMFlift{l}$ maps to the space of $\widetilde{\Sym}^{l}(V)$-valued closed differential forms on $Y_{L}$. 
Here, $\widetilde{\Sym}^{l}(V)$ denotes the local system on $\Grassm$ 
associated to the $l$-th symmetric power of $V$. 
Details are provided in the primary source or in a compact fashion in \cite{BruinierFunkesurj2010} 
and we refer to the lift $\KMFlift{l}$ as the \emph{Kudla--Millson--Funke}\index{Kudla--Millson--Funke!lift}\index{lift!Kudla--Millson--Funke} 
lift of degree $l \in \N_{0}$. 
The latter source also contains the following injectivity statement. 

\begin{thm}[{\cite[Cor.~4.11]{BruinierFunkesurj2010}}]
	Assume that $\dimV > \max\{4, 3 + r \}, \possig > 1, \negsig + l$ even, 
	and that $L$ is even unimodular. 
	Then the theta lift 
	$\KMFlift{l} : \CF_{L,k+l} \to \Dforms^{\negsig} (Y_{L}, \widetilde{\Sym}^{l} (V))$ is injective.
\end{thm}

This approach has recently been generalised by Stein to include the case of maximal lattices 
by solving intricate local integrals for bad primes in the same setting. 

\begin{thm}[{\cite[Cor.~7.8]{Stein2023}}]
	Let $\dimV > \max\{6,2l-2,3+r\}$ and assume that $\negsig + l$ as well as $k = \tfrac{\dimV}{2} + l$ are even. 
	If $L'/L$ is assumed to be anisotropic, then $\KMFlift{l} : \CF_{L,k+l} \to \Dforms^{\negsig} (Y_{L}, \widetilde{\Sym}^{l} (V)) $ is injective. 
\end{thm}

Nevertheless, the case with the most compelling application within this paper remains the $\OG(\possig,2)$ setting 
since it relates to the Borcherds lift -- a connection that is highlighted below. 
In this case, the strongest result available in the literature is attributed to Bruinier 
and proven by refining his approach in \cite{Bruinier2002} through the development of a newform theory for 
vector-valued modular forms in order to relax the assumption of a second hyperbolic split. 

\begin{thm}[{\cite[Thm.~5.3]{BruinierConverse}}]\label{thm:Bruinier_KMinj_conversepaper}
	Let $(L,\qq)$ be an even lattice of signature $(\possig,2)$ such that it splits a hyperbolic plane and a scaled hyperbolic plane $L \simeq K \oplus H(N) \oplus H$ for some $N \in \N$. 
	Then the Kudla--Millson lift $\KMlift : \CF_{L,1 + \possig/2} \to \Hforms^{(1,1)}(Y_{L})$ is injective. 
\end{thm}

It should also be noted that quite recently, Zuffetti and the author have extended this result to the case of general signature, also including the Funke--Millson twist.

\begin{thm}[{\cite[Thm.~6.2]{Metzler2026}}]\label{thm:injindeg1_l}
	Let $L$ be an even indefinite lattice of signature~$(\possig,\negsig)$. 
	\begin{enumerate}[label=\roman*)]
		\item If~$L \simeq K \oplus H(N) \oplus H$ for some even lattice~$K$ and some positive integer~$N$, 
		then the lift~$\KMFlift{l}$ is injective.
		\item Let~$\negsig = 1$.
		If~$L \simeq M \oplus H$ for some positive definite even lattice $M$ and~$M\otimes\Z_p$ splits off a hyperbolic plane for every prime~$p$, then the lift $\KMFlift{l}$ is injective.
		\item Let~$\possig = 1$.
		\label{thm:injindeg1:possig=1_l}
		Then the lift~$\KMFlift{l}$ associated to~$L$ is identically zero.
	\end{enumerate}
\end{thm}

It is noteworthy that this theorem implies all of the injectivity results above. 

The Borcherds and Kudla--Millson lift are closely related -- 
a connection that becomes apparent when following the construction of Borcherds' additive lift presented in 
the celebrated paper \cite{Bruinier2004a}. 
We reduce to the case of $L$ splitting a hyperbolic plane in order to avoid complexity which is not necessary within the scope of the present paper 
and state the following result.

\begin{thm}[{\cite[Thm~4.2 p.~330]{BruinierConverse}}]\label{thm:KMliftinj_Borliftsurj}
	Suppose that $L$ splits a hyperbolic plane, $\possig \geq 2$, and that $\possig$ is greater than the Witt rank of $V$. 
	Then the following are equivalent:
	\begin{enumerate}[label=\roman*)]
		\item 
		The map $\KMlift: \CF_{L,k} \to \Hforms^{(1,1)}(Y_{L})$ is injective. 
		
		\item 
		Every meromorphic modular form $F$ with respect to $\disckernel[L]$ 
		whose divisor is a linear combination of special divisors as in \eqref{eq:Psi_divisor} is 
		(up to a nonzero constant factor) the Borcherds lift $\Psi(z, f)$ of a 
		weakly holomorphic modular form $f \in \MF_{L^{-},2-k}^{!}$ with integral principal part.
	\end{enumerate}
\end{thm}

		\section{Results on theta lifts}\label{sec:KMlift_Injectivity}

Let $(L,\qq)$ be an even $\Z$ lattice of signature $(\possig,2)$ and rank $\dimV$ and set $k \coloneq \dimV/2$. 
In this section we will prove the injectivity of the Kudla--Millson lift introduced in Definition~\ref{de:KMlift}:
\[
	\KMlift\colon \CF_{L,k} \to \HF^{(1,1)}(Y_{L}).
\] 
By Theorem~\ref{thm:KMliftinj_Borliftsurj}, this implies the converse theorem for Borcherds products stated in Theorem~\ref{thm:Intro_Conversetheorem}.
First, note that the lift $\KMlift$ is linear, 
hence it suffices to prove that its kernel is trivial. 
The proof idea is for $f \in \CF_{L,k}$ to integrate the form $\KMlift(f)$ 
over special cycles of $Y_{L}$.  
The nature of these cycles identifies the integrals as special values of a family of \(L\)-series. 
If it is assumed that $f$ is annihilated by the Kudla--Millson lift, then these \(L\)-series necessarily have zeros. 
This will allow us to extract information 
about the Fourier coefficients of the initial form $f$ and to
subsequently conclude that they vanish without exception, i.e.\ $f = 0$, proving the injectivity.
\\
In the following, abbreviate $\Df = L'/L$ and denote the Fourier expansion of $f$ by 
\begin{equation}\label{eq:Fourierexp_f_cuspform_injectivityproof_sec}
	f(\tau) = \sum_{\lambda \in \Df} \sum_{n \in \qd(\lambda) + \Z} a(\lambda,n) \cdot e(n \tau) \efr_{\lambda}. 
\end{equation}

\subsection{Cycle integrals of Kudla--Millson lifts}\label{ssec:Cycle_Integrals_KM_Lifts}

The first major step is to prove that integrals of a target form $\KMlift(f)$ of the Kudla--Millson lift over certain divisors 
equal special values of symmetric square type $L$-functions of $f$. 
The relevant divisors have been presented in Subsection~\ref{ssec:Special_divisor} in \eqref{eq:de:candiv} 
and we briefly recall their shape.
Let $\posvec \in L$ have positive norm, $\Grassm$ denote the Grassmannian of $L$, and $\disckernel[L]$ 
denote the discriminant kernel in $\OG(L)$. If $\Gamma(L)_{\posvec}$ denotes the stabiliser of $\posvec$ in $\disckernel[L]$, then 
\[
	\candiv(\posvec) \coloneq \Gamma(L)_{\posvec} \backslash \Grassm_{\posvec} \to \disckernel[L] \backslash \Grassm = Y_{L} 
\]
defines a (in general relative) cycle of $Y_{L}$, 
where we have set $\Grassm_{\posvec} = \{x \in \Grassm \mid x \subseteq \posvec^{\perp}\}$.
We also require the following definitions from \cite{Metzler2025}. 
\begin{de}\label{de:symmsqLf}
	For $f \in \MF_{L,k}$, a splitting sublattice $M = L_1 \oplus L_2 \leq L$ with definite $L_1 \otimes_{\Z} \R$, 
	an isotropic element $\eta \in (L'\cap L_2')/(L \cap L_2')$, and $t \in \Q^{\times}$ define the formal series 
	\begin{equation}\label{eq:symmsqLfunforsublattice}
		L_{L_1,\eta,t}(f,s) \coloneq \sum_{0 \neq l \in L_1' \cap L'} \frac{a(\overline{l} + \eta, t \qq(l))}{\qq(l)^{s}}.
	\end{equation} 
	Further, for $\lambda \in \Df$ we define the associated \index{Lfunction@$L$-function!vector-valued!symmetric square} \emph{symmetric square} $L$-function and its \emph{specialisation} at some $N \in \N$ as 
	\begin{align}
		L_{(\lambda,t)}(f,s) \coloneq \sum_{n \in \N} \frac{a(n\lambda,n^2t)}{n^s}, \qquad L_{(\lambda,t)}^{N}(f,s) \coloneq \sum_{\substack{n \in \N \\ \gcd(n,N) = 1}} \frac{a(n\lambda,n^2t)}{n^s}.
	\end{align}
	Lastly, for $\posvec \in L'$ with $\qq(\posvec) > 0$, we abbreviate \(L_{\posvec}^{N} \coloneq L_{(\lambda,t)}^{N}\) with $(\lambda,t) = (\overline{\posvec}, \qq(\posvec))$. 
\end{de}

Motivated by Kudla's geometric construction \cite[Sec.~4.3]{Kudla2003},
we consider the double integral 
\begin{align}\label{eq:KMlift_cycle_integral_display_for_proof}
	\int_{\candiv(\posvec)} \int_{\Fd} \langle f(\tau) , \Theta_{L}(\tau,z; \KMform) \rangle v^k \frac{\di u \di v}{v^{2}} \wedge \Kaehlerform^{\possig -2},
\end{align}
where $\Kaehlerform$ is the negative of a Kähler form on $Y_{L}$ as in Definition~\ref{de:Omega_form} and $\Fd = \Mp_{2}(\Z)\backslash \Ha$. 
In order to apply Kudla's observation \cite[Prop.~4.17]{Kudla2003}, several reduction steps are performed. 
Let $L_{\posvec} \coloneq L \cap \posvec^{\perp}$, a lattice of signature $(\possig-1,2)$, 
and define the sublattice $M \coloneq \Z\posvec \oplus L_{\posvec} \leq L$. 
If $\OG(V)_{\posvec}$ denotes the stabiliser of~$\posvec$, we find 
\[
	\disckernel[L_{\posvec}] \simeq \disckernel[M] \cap \OG(V)_{\posvec} \leq \disckernel[L] \cap \OG(V)_{\posvec} = \Gamma(L)_{\posvec}. 
\] 
Here, we have used that the discriminant kernel is inclusion preserving. 
Further, the bound  
$[\Gamma(L)_{\posvec} : \disckernel[L_{\posvec}]] \leq [\disckernel[L] : \disckernel[M]]$ 
necessitates that the multiplicity~$C(\posvec)$ of the covering 
\[
	\candivcover(\posvec) \coloneq \disckernel[L_{\posvec}] \backslash \Grassm_{\posvec} \to \disckernel[L]_{\posvec} \backslash \Grassm = \candiv(\posvec)
\] 
is finite. 
Since the integrand in \eqref{eq:KMlift_cycle_integral_display_for_proof} is invariant under pullbacks from $\disckernel[L]$ (cf.\ Definition~\ref{de:Omega_form} and \eqref{eq:KMform_Oinvariance}), 
we may rewrite the integral as an integral over $\candivcover(\posvec)$, provided we adjust for the prefactor $C(\posvec)^{-1}$. 
 
We now proceed with the computation, postponing a justification of convergence to a later stage. 
The first aim is to isolate the part of the integrand that belongs to the cycle $\candiv(\posvec)$ in order to apply the Siegel--Weil formula. 
\\
By Proposition~\ref{prop:arrowopsrelation} the operator $\uplift{L}{M}$ described in~\eqref{eq:upop_concrete} satisfies
\begin{align*}
	\langle f(\tau) , \Theta_{L}(\tau,z ; \KMform) \rangle = \langle \uplift{L}{M} f(\tau) , \Theta_{M}(\tau,z ; \KMform) \rangle.
\end{align*}
Note that $\Theta_M(\tau,z ;\KMform)$ fulfils the conditions of 
Remark~\ref{rem:Thetafun_decomposes_tp_if_Schwartzform_does} on $\Grassm_{\posvec}$ 
(cf.\ \cite[Thm.~2.1]{Funke2002}). 
As a consequence, this theta function splits as a tensor product and 
by abbreviating $\tilde{f} \coloneq \uplift{L}{M} f$ and writing $\varphi_{\textnormal{KM},L_{\posvec}}$ for the Kudla--Millson Schwartz form on $L_{\posvec} \otimes \R$, we obtain the following identity on $\Grassm_{\posvec}$: 
\begin{align*}
	\langle f(\tau) , \Theta_{L}(\tau,z; \KMform) \rangle = 
	\langle \tilde{f}(\tau) , \Theta_{\Z\posvec}(\tau,z; \varphi_{0,\posvec}) \otimes \Theta_{L_{\posvec}}(\tau,z; \varphi_{\textnormal{KM},L_{\posvec}}) \rangle,
\end{align*} 
where $\varphi_{0,\posvec}(x) = \exp(- 2 \pi \qq(x))$ is the standard Gaussian on the space $\R \posvec$. 
Assuming absolute convergence, we interchange the order of integration in \eqref{eq:KMlift_cycle_integral_display_for_proof} to separate the contribution from \(\candivcover(\posvec)\).  
For that purpose, write $\lambda_1 \in (\Z\posvec)'/\Z\posvec$ as well as 
$\lambda_2 \in (L_{\posvec})'/L_{\posvec}$ 
and find in the notation of Definition~\ref{de:thetafun_Schwartzform_coset} that 
\begin{align*}
	&\int_{\candivcover(\posvec)} \int_{\Fd} \langle f(\tau) , \Theta_{L}(\tau,z; \KMform) \rangle v^k \frac{\di u \di v}{v^{2}} \wedge \Kaehlerform^{\possig -2} \\
	=\,	& \int_{\Fd}  \sum_{\lambda_1,\lambda_2} \tilde{f}_{\lambda_1 \oplus \lambda_2}(\tau) \cdot \overline{\theta_{\Z\posvec,\lambda_1}(\tau; \varphi_{0,\posvec})} 
	\int_{\candivcover(\posvec)} \overline{\theta_{L_\posvec,\lambda_2} (\tau,z ; \varphi_{\textnormal{KM},L_{\posvec}})} \wedge  \Kaehlerform^{\possig -2} v^k \frac{\di u \di v}{v^{2}}. 
\end{align*}
Observe that $\Theta_{\Z\posvec}$ has no $z$-dependence, as the associated Grassmannian is trivial. 
The inner integral over $\candivcover(\posvec)$, 
has already been explicitly related to a standard theta integral by means of 
Corollary~\ref{cor:OCycleint_as_standardtheta} which in turn equals 
a special value of an Eisenstein series by the application of the Siegel--Weil formula 
(cf.\ Theorem~\ref{thm:SiegelWeil}).
Assuming that $\possig - 1$ exceeds the Witt index of the rational quadratic space $V_{\posvec} = L_{\posvec} \otimes \Q$, we obtain   
\begin{align}
	&\int_{\candiv(\posvec)} \int_{\Fd} \langle f(\tau) , \Theta_{L}(\tau,z; \KMform) \rangle v^k \frac{\di u \di v}{v^{2}} \wedge \Kaehlerform^{\possig -2} \nonumber\\
	=\;& C_{L}(\posvec) \cdot \int_{\Fd}  \sum_{\lambda_1,\lambda_2} \tilde{f}_{\lambda_1 \oplus \lambda_2} \overline{ \theta_{\Z\posvec,\lambda_1}(\tau; \varphi_{0,\posvec}) \Eisadel(\tau,s_{0};\tilde{\varphi}_{\textrm{KM},L_\posvec} \otimes \varphi_{\lambda_2}) }	v^k \frac{\di u \di v}{v^{2}} \label{eq:InjHypSplit_computation_1}
\end{align}
for a nonzero explicit constant\footnote{
	In fact, we find that $C_{L}(\posvec)$ is given by the product of the prefactors on the right hand side of \eqref{eq:cor:OCycleint_as_standardtheta} in Corollary~\ref{cor:OCycleint_as_standardtheta} with $C(	\posvec)^{-1}$. 
	} 
$C_{L}(\posvec)$. 
Here, the Eisenstein series is to be understood adelically 
and has been presented in Definition~\ref{de:Eisenstein_adelically_depence_Schwartzform}. 
Further, it is holomorphic at the critical point $s_0 = k - 3/2$ 
while  
$\varphi_{\lambda_2}$ denotes the finite Schwartz form 
represented by the coset of $\lambda_{2}$ 
(cf.\ below~\eqref{eq:intertwiningoperator_finitepart}). 
By \eqref{eq:KMfun_intertwines_to_Phik}, the intertwining operator maps $\tilde{\varphi}_{\textnormal{KM},L_{\posvec}}$ \todo{explain tilde notation}
to the standard section $\Phi_{\infty}^{k-1/2}$ of weight $k-1/2$ in the principal series representation. 
Hence, collecting terms for $\lambda_{1}, \lambda_{2}$ and utilising Definition~\ref{de:Eisensteinseries_VV_adelically}, the inner sum in \eqref{eq:InjHypSplit_computation_1} 
may again be written as a scalar product by using the identification of $\C[\Df]$ with $\Sspace_{\Df}$: 
\[
\sum_{\lambda_1,\lambda_2} \tilde{f}_{\lambda_1 \oplus \lambda_2} \overline{ \theta_{\Z\posvec,\lambda_1}(\tau; \varphi_{0,\posvec}) \Eisadel(\tau,s_{0};\tilde{\varphi}_{\textrm{KM},L_\posvec} \otimes \varphi_{\lambda_2}) }
=
\langle \tilde{f} ,  \Theta_{\Z\posvec}(\tau; \varphi_{0,\posvec}) \otimes \Eisadel_{\hat{L}_{\posvec},k-1/2}(\tau,s_{0}) \rangle.
\]
Here, $\Eisadel_{\hat{L}_{\posvec},k-1/2}$ denotes the adelic Eisenstein series introduced in Definition~\ref{de:Eisensteinseries_VV_adelically}, which corresponds to a classical vector-valued Eisenstein series via Proposition~\ref{prop:Eisenstein_adelic_to_VVMF}. 
In fact, \eqref{eq:prop:Eisenstein_adelic_to_VVMF} reads for $l \in \Z/2$, representing the weight, 
\begin{equation}\label{eq:Eisenstein_adelic_to_classical_finalproof}
	\Eisadel_{\hat{L}_{\posvec},l}(\tau,s_{0}) = E_{L_{\posvec},0,l}\left(\tau,\tfrac{s_{0}+1-l}{2}\right),  
\end{equation}
where $E_{L_{\posvec},0,l}$ is described in Definition~\ref{de:VVEisensteinseriesnonholomorphic}. 
As a consequence, selecting $l = k - 1/2$ yields a new expression for the initial integral in classical terms: 
\begin{align}\label{eq:Peterssonproductassociatedtomffandposvec}
	& \;\int_{\candiv(\posvec)} \int_{\Fd} \langle f(\tau) , \Theta_{L}(\tau,z; \KMform) \rangle v^k \frac{\di u \di v}{v^{2}} \wedge \Kaehlerform^{\possig -2} \nonumber\\
	=& \;C_{L}(\posvec) \cdot \int_{\Fd}  \langle \tilde{f}(\tau) ,  \Theta_{\Z\posvec}(\tau; \varphi_{0,\posvec}) \otimes E_{L_{\posvec},0,k-1/2}(\tau,0) \rangle	v^k \frac{\di u \di v}{v^{2}}.
\end{align}
This expression converges absolutely, as long as the Eisenstein series has no pole in \(s\) (cf.\ Remark~\ref{rem:Eis_locunifbound}), 
and it is suitable for applying 
Proposition~5.6 of \cite{Metzler2025}. 
This yields the following result. 
\begin{prop}\label{prop:RankinSelbergVVMFThetaposdefsplitcrossEisenstein}
		Let the notation be as above. 
		Then we have the identity 
		\begin{equation}\label{eq:RankinSelbergVVMFThetaposdefsplitcrossEisenstein}
			\int_{\Fd} \left\langle \tilde{f} , \Theta_{\Z \posvec} \otimes E_{L_\posvec,0,k - 1/2}(\argdot,s) \right\rangle \cdot \Im^k \di \mu	
			=	\frac{\Gamma(\overline{s}+k-1)}{(4 \pi)^{\overline{s}+k-1}} \cdot L_{\Z \posvec,0,1}(\tilde{f},\overline{s} + k - 1)
		\end{equation}
		of holomorphic functions in $\overline{s}$ for $\Re(s) > \frac{ m_1 - k }{2} + 1 - \boundoptimiser$, 
		where $\boundoptimiser = 1/2$ or $1/4$ depending on whether $\rk(L)$ is even or not.  
	\end{prop}
To ease the notation before continuing, we may select a primitive element $\posvecprim \in (\Z\posvec)' \cap L'$. 
Notice that the Fourier coefficients of cusp forms have a symmetry property enforced by the action of the centre of \(\Mp_{2}(\Z)\), 
yielding that the series vanishes unless $2k \equiv \sig(L) \mod 4$. 
In our case $2k = \possig + 2 \equiv \possig - 2 \mod 4$, which equals the signature. 
Hence, the positive and negative indices of the series contribute symmetrically, yielding 
\begin{align*}
	L_{\Z \posvec, 0, 1}(\tilde{f}, \overline{s} + k -1) 
	&= \frac{2}{\qq(\posvecprim)^{\overline{s}+k-1}} \cdot L_{\posvecprim}(f,2(\overline{s} + k - 1)). 
\end{align*} 
We apply our findings to \eqref{eq:Peterssonproductassociatedtomffandposvec} 
in order to infer the following Theorem. 

\begin{thm}\label{thm:KMlift_cycleint_to_Lvalue}
	Let $(L,\qq)$ be an even lattice of signature $(\possig,2)$, 
	set $k \coloneq 1 + \tfrac{\possig}{2}$, let $f \in \CF_{L,k}$ and select some $\posvec \in L$ of positive norm. 
	If $\possig > 3$, then 
	\begin{equation}
		\int_{\candiv(\posvec)} \KMlift(f) \wedge \Kaehlerform^{\possig -2} = 2C_{L}(\posvec) \frac{\Gamma(k-1)}{(4 \pi\qq(\posvecprim))^{k-1}} \cdot L_{\posvecprim}(f,\possig),
	\end{equation}
	where $\KMlift$ is the Kudla--Millson lift, 
	$\Kaehlerform$ the negative of a Kähler form presented in Definition~\ref{de:Omega_form}, 
	$\candiv(\posvec)$ is the special divisor in \eqref{eq:de:candiv}, 
	$C_{L}(\posvec)$ is a nonzero explicit constant, 
	$\posvecprim$ is some primitive element in $\Q \posvec \cap L'$, and 
	$L_{\posvecprim}(f,s)$ denotes the symmetric square type $L$-series from Definition~\ref{de:symmsqLf}. 
\end{thm}

As a consequence, we conclude that a family of special $L$-values vanishes, 
if the associated cusp form is annihilated by the Kudla--Millson lift.

\begin{cor}\label{cor:KMlift_vanishes_Lvalues_vanish}
	In the notation of Theorem~\ref{thm:KMlift_cycleint_to_Lvalue} we find that for any $f \in \ker(\KMlift)$ 
	and all primitive $\posvecprim \in L'$ of positive norm the following $L$-value vanishes: 
	\begin{equation}\label{eq:cor:KMlift_vanishes_Lvalues_vanish}
		L_{\posvecprim}(f,\possig) \overset{\textnormal{def}}{=} \sum_{ n \in \N } \frac{a(n\overline{\posvecprim} , n^2\qq(\posvecprim))}{n^{\possig}} = 0.
	\end{equation}
\end{cor}

The corollary above will be employed to derive an injectivity result.

\subsection{The case of a hyperbolic split}\label{ssec:Injectivity_hyp_split}

The aim of this subsection is to prove the following theorem. 

\begin{thm}\label{thm:mainhyperbolicsplit}
	Let $(L,\qq)$ be an even lattice of signature $(\possig,2)$ with $\possig > 3$. 
	Assume that $L$ splits a hyperbolic plane. 
	Then the Kudla--Millson lift ${\KMlift}\colon \CF_{L,1+\possig/2} \to \mathcal{H}^{(1,1)}(Y_{L})$ is injective.
\end{thm}

Recall that the Kudla--Millson lift associates to 
an elliptic cusp form $f \in \CF_{L,k}$ a form 
\[
	\KMlift (f) = \int_{\Mp_{2}(\Z) \backslash \Ha} v^k \langle f(\tau) , \Theta(\tau,z; \varphi_{\textnormal{KM}}) \rangle \frac{\di u \di v}{v^{2}} \in \mathcal{H}^{(1,1)}(Y_{L}).
\]
\todocre{improve part above}

We now outline the strategy of the proof. 
By Corollary~\ref{cor:KMlift_vanishes_Lvalues_vanish}, the vanishing of $\KMlift(f)$ for some 
$f \in \CF_{L,k}$ implies the vanishing of special $L$-values associated 
to primitive $\posvecprim \in L'$ of positive norm. 
These \(L\)-values are constructed by summing over Fourier coefficients supported along rational lines in $V$. 
The key observation is that if \(L\) splits a hyperbolic plane, then the primitivity assumption may be dropped.
With that advantage, the series $L_{\posvecprim}(f,\possig)$ may be decomposed into vanishing subseries, 
supported on arbitrary multiples of \(\posvecprim\), which suffices to force the vanishing of the Fourier coefficients themselves.

In order to deduce the vanishing of \(f\) from the vanishing of the \(L\)-values in \eqref{eq:cor:KMlift_vanishes_Lvalues_vanish} two issues must be addressed. 
First, the \(L\)-functions involve only coefficients of the form $a(\overline{\posvec}, \qq(\posvec))$ for vectors $\posvec \in L'$ of positive norm, 
leaving open the behaviour of coefficients with other indices. 
Second, the resulting linear relations span an infinite dimensional system and are not obviously sufficient to force \(f = 0\). 
While the second issue can, in principle, be resolved using Hecke theory, the first is solely a property of the lattice \(L\). 
The key insight is that under a hyperbolic split, all indices can be realised as \((\overline{\posvec}, q(\posvec))\) for suitable \( \posvec \in L'\).

\begin{lemma}\label{lem:Hypsplit_Indices_primitivelyrep}
	Suppose that $L \simeq K \oplus H$ splits a hyperbolic plane. Then every pair $(\lambda,n)$, where \(\lambda \in \Df\) and \(n \in \qd(\lambda) + \Z\) 
	arises as $(\overline{\posvecprim},\qq(\posvecprim))$ 
	for some primitive element $\posvecprim \in L'$. 
\end{lemma}

\begin{proof}
	Let $(\lambda,n)$ be an admissible index for a Fourier coefficient.
	Choose a representative $\ell \in L'$ with $\overline{\ell} = \lambda$, 
	meaning $\qq(\ell) \equiv n \mod \Z$ 
	and note that the class $\overline{\ell} \in \Df$ is invariant under translating $\ell$ by elements in $H$. 
	Since $L \simeq K \oplus H$ splits a hyperbolic plane $H = \langle e_1, e_2\rangle_\Z$, 
	where $e_1,e_2$ denote the standard generators with $\qq(e_1) = \qq(e_2) = 0$ and $\bilf(e_1,e_2)=1$, we find that \( \qq(\ell) = \qq(\ell\vert_{K}) + \qq(\ell\vert_{H}) \)
	and for all $m_1, m_2 \in \Z$ that \( \qq(m_1 e_1 + m_2 e_2) = m_1 m_2 \)
	so that the primitive vector 
	\[
		\ell_\lambda \coloneq \ell\vert_{K} \oplus \left[(n - \qq(\ell\vert_{K})) e_1 + e_2\right] 
	\]
	fulfils the desired property $(\overline{\ell_\lambda}, \qq(\ell_\lambda)) = (\lambda,n)$. 
\end{proof}

We note that the $L$-series $L_{\posvec}(f,s)$ depends solely on the pair $(\overline{\posvec},\qq(\posvec))$, 
so that in conjunction with Corollary~\ref{cor:KMlift_vanishes_Lvalues_vanish} we immediately conclude the following assertion. 

\begin{cor}\label{cor:Lvalue_positivevec_vanish_for_hypsplit}
	Let $(L,\qq)$ be an even $\Z$ lattice of signature $(\possig,2)$ and assume
	it splits a hyperbolic plane. 
	Set $k \coloneq 1 + \possig/2$ and suppose that $f \in \CF_{L,k}$ is annihilated by $\KMlift$. 
	If $\possig > 3$, then for any $\ell \in L'$ of positive norm the following special $L$-value vanishes:
	\[
	L_{\posvec}(f,\possig) \overset{\textnormal{def}}{=} \sum_{n \in \N}  \frac{a( n \overline{\posvec}, \qq(n\posvec))}{n^{\possig}} = 0.
	\]
\end{cor}

As we will see, this also eliminates the second problem in proving the main theorem without resorting to Hecke theory by utilising an inclusion--exclusion trick. 

\begin{lemma}\label{lem:LvalueequalspartialLvalueifhypsplit}
	Assume $L$ splits a hyperbolic plane and $f \in \CF_{L,k}$ for $\possig > 3$ is annihilated by $\KMlift$.  
	Let $N\in \N$ be a natural number and $\posvec \in L'$ be of positive norm. 
	Then
	\[
	L_{\posvec}^N(f,\possig) \overset{\textnormal{def}}{=} \sum_{\substack{n \in \N \\ \gcd(n , N) = 1}}  \frac{a( n \overline{\posvec}, \qq(n\posvec))}{n^{\possig}} = 0.
	\]
\end{lemma}

\begin{proof}
	Note that it suffices to consider square free $N$. 
	We will prove the assertion by induction on the number of distinct prime divisors of $N$. 
	The case $N = 1$ is identical to Corollary~\ref{cor:Lvalue_positivevec_vanish_for_hypsplit}. 
	Suppose the statement is true for some $N \in \N$. 
	Select a prime number~$p \nmid N$ and note that by the assumption $\possig > 3$ 
	absolute convergence guarantees 
	\begin{align*}
	\sum_{\substack{n \in \N \\ \gcd(n , pN) = 1}}  \frac{a( n \overline{\posvec}, \qq(n\posvec))}{n^{\possig}}\;
	&=
	\sum_{\substack{n \in \N \\ \gcd(n , N) = 1}}  \frac{a( n \overline{\posvec}, \qq(n\posvec))}{n^{\possig}}
	- \frac{1}{p^{\possig}} \sum_{\substack{n \in \N \\ \gcd(n , N) = 1}}  \frac{a( n p\overline{\posvec}, \qq(np\posvec))}{n^{\possig}} \\
	&= L_{\posvec}^{N}(f,\possig) - \frac{1}{p^{\possig}} \, L_{p\posvec}^{N}(f,\possig).
	\end{align*}
	Now the right hand side vanishes by the induction hypothesis, completing the proof. 
\end{proof}

The lemma above is the last ingredient for the proof of the main theorem.

\begin{proof}[{\textit{Proof of Theorem~\ref{thm:mainhyperbolicsplit}}}]
	Suppose there was a non-trivial form $f \in \ker(\KMlift)$ with Fourier expansion 
	\[
	f(\tau) = \sum_{\lambda \in \Df} \sum_{n \in \qd(\lambda) + \Z} a(\lambda,n) \cdot e(n\tau) \efr_{\lambda}.
	\] 
	Select a pair $(\lambda,n)$ for which $a(\lambda,n) \neq 0$. 
	By Lemma~\ref{lem:Hypsplit_Indices_primitivelyrep}, the hyperbolic split guarantees the existence of some $\posvec \in L'$ of positive norm with $(\overline{\posvec}, \qq(\posvec)) = (\lambda,n)$. 
	Further, the assumption of~$\possig > 3$ implies that $L_{\posvec}(f,\possig) = \sum_{d \in \N} a(d \overline{\posvec},d^2\qq(\posvec))/d^{\possig}$ converges absolutely. 
	Hence, there must be some number $M \in \N$ such that 
	\[
	\sum_{\substack{d \in \N \\ d > M}}^\infty \frac{\lvert a(d \overline{\posvec},d^2 \qq(\posvec))\rvert}{d^\possig} < \lvert a(\overline{\posvec},\qq(\posvec)) \rvert.
	\]
	Consequently, selecting $N \coloneq M!$ results in 
	\begin{equation}\label{eq:MainProof_1}
		L_{\posvec}^{N}(f,\possig) 		= a(\overline{\posvec},\qq(\posvec)) + \sum_{\substack{1 < d \in \N \\ \gcd(d,N) = 1}} \frac{a(d \overline{\posvec}, d^2 \qq(\posvec))}{d^{\possig}} \neq 0.
	\end{equation} 
	However, $f$ was assumed to be annihilated by $\KMlift$, implying by Lemma~\ref{lem:LvalueequalspartialLvalueifhypsplit} that the left hand side of \eqref{eq:MainProof_1} vanishes 
	-- a contradiction!
						As a consequence, the assumption of a non-trivial element in the kernel of $\KMlift$ must be incorrect i.e.\ 
	the Kudla--Millson lift is injective.  
\end{proof}

An application of Theorem~\ref{thm:KMliftinj_Borliftsurj} yields a converse theorem for Borcherds products (also compare Theorem~\ref{thm:converse_theorem_Borcherds_products_possiblesublattice}). 

\begin{thm}\label{thm:converse_theorem_Borcherds_products}
	Assume that $L \simeq K \oplus H$ for some lattice $K$ of signature $(\possig-1,1)$ with $\possig > 3$. 
	Then every meromorphic modular form $F$ for $\disckernel[L]$ whose divisor is a linear combination of special divisors $Z(\mu,n)$ as in \eqref{eq:de:Heegner_divisor} is (up to a nonzero constant factor) the Borcherds lift $\Psi(z,f)$ 
	of some weakly holomorphic modular form $f \in \MF_{L^-,1-\possig/2}^{!}$. 
\end{thm}

\begin{rem}
	The assumption of a hyperbolic split is essential for Theorem~\ref{thm:mainhyperbolicsplit} and may not be omitted. 
	This has been demonstrated in \cite[Sec.~6.1 p.~333]{BruinierConverse}, where counterexamples are constructed. 
	Moreover, the theorem does not extend to the case $\possig = 1$; see\ Theorem~\ref{thm:injindeg1_l}~\ref{thm:injindeg1:possig=1_l}. 
	The case of $\possig = 3$ is discussed in the following, while the case of $\possig = 2$ 
	remains inaccessible for a single hyperbolic split with the current method. 
\end{rem}

\subsection{The case of $\possig = 3$}

The case of $\possig = 3$, which was previously excluded for the sake of convenience, 
may still be treated by the same procedure. 
A careful analysis of the proof reveals that the majority of arguments remains intact 
if vanishing of the $L$-value $L_{\posvec}(f,\possig)$ is replaced by 
convergence $L_{\posvec}(f,s) \to 0$ for $s \to \possig$. 
This approach is universally applicable for Witt rank $1$ of the lattice $L$ 
which is precisely the case that remains unproven in the current literature. 
The crucial additional element for such an advancement is Corollary~5.4 from \cite{Metzler2025}, 
stating that \(L_{\posvec}^{N}(f,s)\) converges absolutely for \(\Re(s) > k + 1 - 2\sigma\) with 
\(\sigma = 5/16\) or \(1/4\), based upon whether \(\level(L) \mid N^{\infty}\) or not. 

Indeed, the sole issue with the case $\possig = 3$ in the proof of Theorem~\ref{thm:mainhyperbolicsplit} 
is convergence. 
More explicitly, in order to apply the Siegel--Weil formula to insert the Eisenstein series~$\Eisadel$ in \eqref{eq:InjHypSplit_computation_1} and obtain holomorphicity at $s_0 = k - 1/2$, Weil's convergence criterion is applied.
This criterion fails, unless $L \cap \posvec^{\perp}$ has Witt rank at most $1$ 
which is guaranteed if~$L$ itself is assumed to have Witt rank $1$. 
Continuing the proof under this assumption, yields that \eqref{eq:Peterssonproductassociatedtomffandposvec}, 
when inserting a general argument \(s\) into the Eisenstein series involved, equals  
\[
	2 \frac{\Gamma(\overline{s}+k-1)}{(\qq(\posvecprim)4 \pi)^{\overline{s}+k-1}} \cdot L_{\posvecprim}(f,2(\overline{s} + k - 1))
\]
for all $s \in \C$ in a right half plane, where the $L$-series converges absolutely and is holomorphic in $\overline{s}$. 
This is the case for $\Re(s) > 0$. 
However, the expression~\eqref{eq:Peterssonproductassociatedtomffandposvec} is holomorphic at $s = 0$ by Weil's convergence criterion. 
Hence, the right hand side of \eqref{eq:RankinSelbergVVMFThetaposdefsplitcrossEisenstein} can be continued to $s = 0$. 
As a consequence, we obtain the following version of Corollary~\ref{cor:KMlift_vanishes_Lvalues_vanish}.

\begin{cor}\label{cor:KMlift_cycleint_to_Lvalue_possig3}
	Let $(L,\qq)$ be an even lattice of signature $(3,2)$, 
	set $k \coloneq \tfrac{5}{2}$, let $f \in \CF_{L,k}$ and select some $\posvec \in L$ of positive norm 
	such that $L \cap \posvec^{\perp}$ has Witt rank at most $1$. 
	Then 
	\begin{equation}
		\int_{\candiv(\posvec)} \KMlift(f) \wedge \Kaehlerform^{1} = 2C_{L}(\posvec) \frac{\Gamma(k-1)}{(4 \pi\qq(\posvecprim))^{k-1}} \cdot \lim_{\substack{s \to \possig \\ \Re(s) > \possig}} L_{\posvecprim}(f,s),
	\end{equation}
	for some primitive element $\posvecprim \in \Q \posvec \cap L'$. 
	For the rest of the notation, compare Theorem~\ref{thm:KMlift_cycleint_to_Lvalue}. 
\end{cor}

\begin{cor}\label{cor:KMlift_vanishes_Lvalues_vanish_possig3}
	In the notation of Corollary~\ref{cor:KMlift_cycleint_to_Lvalue_possig3}, 
	assume $f \in \ker(\KMlift)$. 
	Then 
	\begin{equation}
		\lim_{\substack{s \to \possig \\ \Re(s) > \possig}} L_{\posvecprim}(f,s) \overset{\textnormal{def}}{=} \lim_{\substack{s \to \possig \\ \Re(s) > \possig}} \sum_{ n \in \N } \frac{a(n\overline{\posvecprim} , n^2\qq(\posvecprim))}{n^{s}} = 0,
	\end{equation}
	where $a(\lambda,n)$ denote the Fourier coefficients of $f$ as in \eqref{eq:Fourierexp_f_cuspform_injectivityproof_sec}. 
\end{cor}

The same line of reasoning that was used to prove Corollary~\ref{cor:Lvalue_positivevec_vanish_for_hypsplit} and Lemma~\ref{lem:LvalueequalspartialLvalueifhypsplit} can be applied to prove their respective analogues 
which are subsumed in the following Corollary. 

\begin{cor}\label{cor:Lvalue_equals_partialLvalue_hypsplit_possig3}
	Let $(L,\qq)$ be an even lattice of signature $(3,2)$ with Witt rank $1$ and assume that it splits a hyperbolic plane. 
	Set \(k=5/2\) and assume that $f \in \CF_{L,k}$ is annihilated by $\KMlift$.  
	Then for any $\posvec \in L'$ be of positive norm and $N\in \N$, we have  
	\[
	\lim_{\substack{s \to \possig \\ \Re(s) > \possig}} L_{\posvec}(f,s) = \lim_{\substack{s \to \possig \\ \Re(s) > \possig}} L_{\posvec}^N(f,s) = 0.
	\]
\end{cor}

\begin{prop}\label{prop:mainhyperbolicsplit_possig3}
	Let $(L,\qq)$ be an even lattice of Witt rank $1$ and signature $(3,2)$, 
	and assume that $L$ splits a hyperbolic plane.  
	Then the lift $\KMlift\colon \CF_{L,5/2} \to \mathcal{H}^{(1,1)}(Y_{L})$ is injective.
\end{prop}

\begin{proof}
	The argument follows that of Theorem~\ref{thm:mainhyperbolicsplit}, with minor modifications to account for 
	a limit process at the critical point $s = 3$ of the $L$-series involved. 
	\\
	Suppose $f \in \ker(\KMlift)$ was non-trivial 
	and choose a non-vanishing Fourier coefficient \(a(\lambda,n) \neq 0 \). 
	By Lemma~\ref{lem:Hypsplit_Indices_primitivelyrep}, 
	the hyperbolic split guarantees the existence of some element $\posvec \in L'$ of positive norm with $(\overline{\posvec}, \qq(\posvec)) = (\lambda,n)$. 
	Further, $\possig = 3$ is chosen such that $L_{\posvec}(f,s) = \sum_{d \in \N} a(d \overline{\posvec},d^2\qq(\posvec))/d^{s}$ converges absolutely to the right of $s = \possig$. 
	In addition, for $N \in \N$ with $\level(L) \mid N^{\infty}$, the series 
	\[
		L_{\posvec}^{N}(f,s) = \sum_{\substack{d \in \N \\ \gcd(d,N) = 1}} 
		\frac{a(d \overline{\posvec},d^2\qq(\posvec))}{d^{s}} 
	\]
	converges absolutely at $s = \possig$. 
	Hence, there must be some number $M \in \N$ such that 
	\[
		\sum_{\substack{M < d \in \N \\ \gcd(d,N) = 1}}^\infty \frac{\lvert a(d \overline{\posvec},d^2 \qq(\posvec))\rvert}{d^\possig} < \lvert a(\overline{\posvec},\qq(\posvec)) \rvert.
	\]
	Consequently, selecting $N \coloneq (M \cdot \level(L))!$ results in 
	\begin{equation}\label{eq:MainProof_1_possig3}
		L_{\posvec}^{N}(f,\possig) 		= a(\overline{\posvec},\qq(\posvec)) + \sum_{\substack{1 < d \in \N \\ \gcd(d,N) = 1}} \frac{a(d \overline{\posvec}, d^2 \qq(\posvec))}{d^{\possig}} \neq 0.
	\end{equation} 
	However, the form $f$ was assumed to be annihilated by $\KMlift$, implying by Corollary~\ref{cor:Lvalue_equals_partialLvalue_hypsplit_possig3} that 
	\[
		L_{\posvec}^N(f,\possig) = \lim_{\substack{s \to \possig \\ \Re(s) > \possig}} L_{\posvec}(f,s) = 0,
	\]
	contradicting \eqref{eq:MainProof_1_possig3}. 
	As a consequence, the assumption of a non-trivial element in the kernel of $\KMlift$ must be incorrect, i.e.\ 
	the Kudla--Millson lift is injective. 
\end{proof}

The preceding result can be employed to augment the converse theorem for Borcherds products. 
Recall that Theorem~\ref{thm:converse_theorem_Borcherds_products} already establishes the case $\possig > 3$ under a single hyperbolic split. 
When $\possig = 3$ and \(L\) admits a hyperbolic split, it necessarily has Witt rank~$1$ or~$2$. 
The case of Witt rank~$2$ is treated in \cite[Thm.~1.2]{BruinierConverse} with an additional technical condition.  
The case of Witt rank $1$, which remained previously unresolved, now follows from Proposition~\ref{prop:mainhyperbolicsplit_possig3} 
in conjunction with Theorem~\ref{thm:KMliftinj_Borliftsurj}.

\begin{thm}\label{thm:converse_theorem_Borcherds_products_possiblesublattice}
	Assume $\possig > 2$ and that $L \simeq K \oplus H$ splits a hyperbolic plane. 
	Then there is a sublattice $K_0 \leq K$	such that every meromorphic modular form $F$ with respect to $\disckernel[L]$ whose divisor is a linear combination of special divisors $Z(\mu,n)$ as in \eqref{eq:de:Heegner_divisor} is (up to a nonzero constant factor) the Borcherds lift $\Psi(z,f)$ 
	of a weakly holomorphic modular form $f \in \MF_{K_0^-,1-\possig/2}^{!}$. 
\end{thm}

\begin{rem}
	The only scenario in which passing to a proper sublattice $K_0 < K$ in Theorem~\ref{thm:converse_theorem_Borcherds_products_possiblesublattice} 
	might be required is when $\possig = 3$ and $L$ has Witt rank $2$. 
	However, we anticipate the method of proof employed in the current paper 
	to also work in this case by bypassing Weil's convergence criterion, 
	rendering the transition to a sublattice $K_0$ superfluous in each instance.  
\end{rem}

		\printbibliography
\end{document}